\documentclass[a4paper,12pt]{article}
\usepackage{amsmath, amssymb, latexsym, amscd, amsthm,amsfonts,amstext}
\usepackage[mathscr]{eucal}
\usepackage{graphicx}
\usepackage{subfig}

\numberwithin{equation}{section}
\input{tcilatex}
\begin{document}

\title{A H\"{o}lder Stability Estimate for a 3D Coefficient Inverse Problem
for a Hyperbolic Equation With a Plane Wave}
\author{ Michael V. Klibanov$^{\ast \text{ }}$and Vladimir G. Romanov$^{\ast
\ast }$ \and $^{\ast }$Department of Mathematics and Statistics, University
of North \and Carolina at Charlotte, Charlotte, NC 28223, USA \and $^{\ast
\ast }$Sobolev Institute of Mathematics, Novosibirsk 630090, \and Russian
Federation \and \texttt{mklibanv@uncc.edu}, \texttt{romanov{@}math.nsc.ru} }
\date{}
\maketitle

\begin{abstract}
A 3D coefficient inverse problem for a hyperbolic equation with
non-overdetermined data is considered. The forward problem is the Cauchy
problems with the initial condition the delta function concentrated at a
single plane (i.e. the plane wave). A certain associated operator is written
in finite differences with respect to two out of three spatial variables,
i.e. \textquotedblleft partial finite differences". The grid step size is
bounded from the below by a fixed number. A Carleman estimate is applied to
obtain, for the first time, a H\"{o}lder stability estimates for this
problem. Another new result is an estimate from the below of the amplitude
of the first term of the expansion of the solution of the forward problem
near the characteristic wedge.
\end{abstract}

\textbf{Key Words}: coefficient inverse problem, hyperbolic equation,
geodesic lines, Carleman estimate, H\"{o}lder stability estimate.

\textbf{2010 Mathematics Subject Classification:} 35R30.

\section{Introduction}

\label{sec:1}

We consider a forward Cauchy problem for a wave-like PDE in $\mathbb{R}^{3}.$
In this problem, the source is the $\delta -$function concentrated at a
plane, which models the incident plane wave. For this problems, we consider
a Coefficient Inverse Problem (CIP). Applications of this CIP is discussed
in section 2. We obtain a H\"{o}lder stability estimate for our CIP.
Uniqueness theorem follows immediately from this estimate. 

We assume that a certain $2\times 2$ system of non local PDEs associated
with the original Cauchy problem is written in partial finite differences,
i.e. finite differences with respect to two out of three spatial variables.
The derivatives with respect to both the third spatial variable and time are
written in the conventional manner. The grid step size of finite differences
is assumed to be bounded from the below by a certain small positive
constant. We point out that such a bound always takes place in computations.
Thus, this assumption has an applied meaning.

The data for our CIP are non-overdetermined ones, i.e. the number $m=3$ of
free variables in the data equals the number $n=3$ of free variables in the
unknown coefficient. It is well known that uniqueness and stability results
for multidimensional CIPs with non-overdetermined data are currently proven
only by the method of \cite{BukhKlib}, which is based on Carleman estimates,
see, e.g. \cite{BK,BY,KT,KL,Yam} and references cited therein for some
samples of publications, which use this method.

Applications of the idea of \cite{BukhKlib} to the developments of the
so-called \textquotedblleft convexification" globally convergent numerical
method for CIPs can be found in, e.g. \cite{Khoa1,Khoa2,KLZ,Klib1d,KL}. The
originating publications about the convexification are \cite{KlibIous,Klib97}%
. Numerical studies by the convexification method of a similar CIP for the
same hyperbolic PDE as the one of this paper and with the point source in
the initial condition can be found in \cite{KLZ} and \cite[Chapter 8]{KL}.
The second generation of the convexification method was developed in \cite%
{Baud1,Baud2,Le}.

However, in terms of uniqueness theorems and stability estimates, the method
of \cite{BukhKlib} works only under the assumption that one of initial
conditions of a corresponding PDE does not equal zero in the entire domain
of interest. The only exceptions are two recent works of Rakesh and Salo 
\cite{RS1,RS2} as well as some follows up publications of these authors. In 
\cite{RS1} a stability result was proven for a CIP for the PDE $%
u_{tt}=\Delta u+q\left( \mathbf{x}\right) u,\mathbf{x}\in \mathbb{R}^{3}$
with two initializing plane waves and with non-overdetermined data. In \cite%
{RS2} uniqueness is proven for a CIP with non-overdetermined data for the
analog of this equation in the frequency domain.

In this paper, we consider the case of the hyperbolic PDE with a non
constant unknown coefficient in the principal part of the hyperbolic
operator with a single incident plane wave. Uniqueness and stability results
were not proven for this CIP in the past. Our case is more complicated than
the one in \cite{RS1}. This is basically because the geodesic lines in our
case are curves rather than straight lines of the above PDE. The above
mentioned assumption of partial finite differences for boundary value
problems for $2\times 2$ systems of non-local PDEs is imposed to avoid the
assumption of \cite{BukhKlib} of the non-vanishing initial condition.

To prove our target result, we use the expansion of the solution of our
forward problem near the characteristic wedge. Such an expansion is known
only for the case of the point source \cite[Theorem 4.1]{Rom1}, \cite[Lemma
2.2.1]{Rom2}. Then we combine this expansion with the Carleman estimate for
the above mentioned $2\times 2$ systems of non-local PDEs.

There are two new results of this paper:

\begin{enumerate}
\item An estimate from the below of the amplitude of the first term of the
expansion of the solution of the forward problem near the characteristic
wedge (Theorem 3.1).

\item A H\"{o}lder stability estimate for the above mentioned CIP (Theorem
6.1).
\end{enumerate}

\textbf{Remark 1.1.} \emph{It seems to be convenient to formulate our main
theorems in the next section 2. However, prior to these formulations, we
need to apply some transformations to the solution of our Forward Problem.
Thus, we postpone those formulations until sections 3, 5 and 6. }

The paper is organized as follows. In section 2 we pose forward and inverse
problems and describe some applications of our CIP. In section 3, we derive
the structure of the solution for the forward problem. In section 4, we
change variables and reduce our CIP to a boundary value problem for a $%
2\times 2$ system of non-local PDEs and reformulate in partial finite
differences a boundary value problem, which was derived in this section. In
section 5 we formulate and prove \noindent two Carleman estimates. In
section 6 we prove a stability estimate for the CIP.

\section{Statements of Forward and Inverse Problems}

\label{sec:2}

Below $\mathbf{x}=\left( x,y,z\right) \in \mathbb{R}^{3}.$ Denote 
\begin{equation*}
\mathbb{R}_{+}^{3}=\left\{ \left( x,y,z\right) :z>0\right\} ,\text{ }\mathbb{%
R}_{-}^{3}=\left\{ \left( x,y,z\right) :z<0\right\} ,\text{ }\mathbb{R}%
_{T}^{4}=\mathbb{R}^{3}\times \left( 0,T\right) .
\end{equation*}%
Let $X,T>0$ be two numbers. We define our domain of interest $\Omega \subset 
\mathbb{R}^{3}$ and related surfaces as:%
\begin{equation}
\Omega =\left\{ \mathbf{x}:-X<x,y<X,~z\in \left( 0,1\right) \right\} ,
\label{2.1}
\end{equation}%
\begin{equation}
\Gamma =\left\{ \mathbf{x=}\left( x,y,z\right) :-X<x,y<X,z=0\right\} \subset
\partial \Omega ,~\Gamma _{T}=\Gamma \times \left( 0,T\right) ,  \label{2.6}
\end{equation}%
\begin{equation}
\Gamma ^{\prime }=\left\{ \mathbf{x=}\left( x,y,z\right)
:-X<x,y<X,z=1\right\} \subset \partial \Omega ,~\Gamma _{T}^{\prime }=\Gamma
^{\prime }\times \left( 0,T\right) ,  \label{2.60}
\end{equation}%
\begin{equation}
\Theta =\left\{ x,y=\pm X,z\in \left( 0,1\right) \right\} \subset \partial
\Omega ,\text{ }\Theta _{T}=\Theta \times \left( 0,T\right) .  \label{2.70}
\end{equation}%
We assume that the function $n\left( \mathbf{x}\right) $ be defined in $%
\mathbb{R}^{3}$ and 
\begin{equation}
n\left( \mathbf{x}\right) \in C^{16}\left( \mathbb{R}^{3}\right) .
\label{3.1}
\end{equation}%
Even though the smoothness requirement (\ref{3.1}) seems to be excessive, we
are unaware how to decrease the smoothness. Indeed, we pose the Forward
Problem in section 2. Next, we prove in section 3 a certain representation
of the solution of this problem. This representation uses (\ref{3.1}). In
addition, it is worth to point out that the issue of the minimal smoothness
is traditionally a minor concern in the field of Coefficient Inverse
Problems, see, e.g. \cite{Nov1,Nov2}, \cite[Theorem 4.1]{Rom1}.

Consider two arbitrary numbers $n_{0},n_{00}$ such that $1<n_{0}<n_{00}.$ We
assume everywhere below that 
\begin{equation}
1\leq n(\mathbf{x})\leq n_{0},~\mathbf{x}\in \mathbb{R}^{3},\Vert n\Vert
_{C^{2}(\mathbb{R}^{3})}\leq n_{00},  \label{2.2}
\end{equation}%
\begin{equation}
n\left( \mathbf{x}\right) =1\text{ in }\mathbb{R}_{-}^{3}\cup \{~|x|\geq
X,~|y|\geq X\}.  \label{2.3}
\end{equation}

In the case of an electric wave field propagation, $n\left( \mathbf{x}%
\right) $ is the refractive index, $n^{2}\left( \mathbf{x}\right) $ is the
spatially distributed dielectric constant of the medium, and the function $%
u\left( \mathbf{x},t\right) $, which we introduce below, is a component of
that electric field. Therefore, the CIP, which we study here, has an
application in the problem of the determination of the dielectric constant
of the medium using the data of the scattered electric wave. These data are
measured at a part $\left( \Gamma \cup \Theta \right) \subset \partial
\Omega $ of the boundary $\partial \Omega $ of the domain $\Omega $. It was
shown numerically on microwave experimental data in \cite{Khoa1,Khoa2,KL}
and analytically in \cite{RK} that in some applications, one wave-like PDE
can govern the process of electric waves propagations equally well with the
full system of Maxwell's equations. In fact, this is also claimed
heuristically in the classic textbook of Born and Wolf \cite[pages 695,696]%
{BW}. Another area of applications of the CIP of this paper is acoustics, in
which case the function $1/n\left( \mathbf{x}\right) $ is the speed of sound
waves propagation in the medium.

\textbf{Forward Problem}. \emph{Solve the following Cauchy problem:}%
\begin{equation}
n^{2}\left( \mathbf{x}\right) u_{tt}-\Delta u=\delta \left( z\right) \delta
\left( t\right) ,\text{ }\left( \mathbf{x},t\right) \in \mathbb{R}^{4},
\label{2.12}
\end{equation}%
\begin{equation}
u\mathbf{\mid }_{t<0}=0.  \label{2.13}
\end{equation}

Since by (\ref{2.3}) $n\left( \mathbf{x}\right) =1$ in $\mathbb{R}_{-}^{3},$
then the incident wave in (\ref{2.12}), (\ref{2.13}) is the plane wave,
which is incident at the plane $\Sigma ,$%
\begin{equation}
\Sigma =\left\{ \mathbf{x}\in \mathbb{R}^{3}\,:z=0\right\}  \label{3.200}
\end{equation}%
and propagates in the direction parallel to the $z-$axis. The function $%
u\left( \mathbf{x},t\right) $ in the CIP is the solution of the Forward
Problem.

\textbf{Coefficient Inverse Problem} (CIP). \emph{Suppose that the following
functions}

\emph{\ }$f_{0},f_{1},f_{2}$\emph{\ are given:}%
\begin{equation}
u\mid _{\Gamma _{T}}=f_{0}\left( \mathbf{x},t\right) ,u_{z}\mid _{\Gamma
_{T}}=f_{1}\left( \mathbf{x},t\right) ,u\mid _{\Theta _{T}}=f_{2}\left( 
\mathbf{x},t\right) ,  \label{2.14}
\end{equation}%
\emph{where }$u\left( \mathbf{x},t\right) $ \emph{is the solution of Forward
Problem (\ref{2.12}), (\ref{2.13}). Determine the function }$n\left( \mathbf{%
x}\right) $ \emph{for} $\mathbf{x}\in \Omega ,$ \emph{assuming that this
function satisfies conditions (\ref{2.2}), (\ref{2.3}).}

It follows from (\ref{2.1}), (\ref{2.6}) that the data at $\Gamma _{T}$ in (%
\ref{2.14}) are the backscattering data.

\textbf{Remark 2.1.} \emph{Let} $S_{T}=\partial \Omega \times \left(
0,T\right) .$ \emph{Note that we do not need the data at }$S_{T}\diagdown
\left( \Gamma _{T}\cup \Theta _{T}\right) $ \emph{in our CIP,\ i.e. we do
not need the data on the transmitted side of the rectangular prism }$\Omega
. $

\section{The Structure of the Solution of the Forward Problem}

\label{sec:3}

\subsection{Geodesic lines}

\label{sec:3.1}

Let $\tau (\mathbf{x})$ be the travel time, which the above plane wave needs
to travel from the plane (\ref{3.200}) to the point $\mathbf{x}\in \mathbb{R}%
_{+}^{3}.$ The function $\tau (\mathbf{x})$ satisfies the eikonal equation
and it equals zero at the plane $\Sigma $, i.e. 
\begin{equation}
\left\vert \nabla \tau (\mathbf{x})\right\vert ^{2}=n^{2}\left( \mathbf{x}%
\right) ,~\mathbf{x}\in \mathbb{R}_{+}^{3};\text{ }\tau |_{z=0}=0.
\label{3.21}
\end{equation}%
We obtain from (\ref{3.21})%
\begin{equation*}
\tau (\mathbf{x})=\pm \dint\limits_{0}^{z}\sqrt{n^{2}\left( x,y,s\right)
-\left( \tau _{x}^{2}+\tau _{y}^{2}\right) \left( x,y,s\right) }\,ds.
\end{equation*}%
However, since Physics tells us that $\tau (\mathbf{x})>0,$ then the correct
formula for $\tau (\mathbf{x})$ is 
\begin{equation}
\tau (\mathbf{x})=\dint\limits_{0}^{z}\sqrt{n^{2}\left( x,y,s\right) -\left(
\tau _{x}^{2}+\tau _{y}^{2}\right) \left( x,y,s\right) }\,ds>0,\text{ }%
\mathbf{x}\in \mathbb{R}_{+}^{3}.  \label{3.2100}
\end{equation}%
In particular, 
\begin{equation}
\partial _{z}\tau (\mathbf{x})>0,\text{ }\mathbf{x}\in \mathbb{R}_{+}^{3}.
\label{3.210}
\end{equation}

The function $\tau(\mathbf{x})$ is generated by the conformal Riemannian
metric with the element of its length 
\begin{equation}
d\tau=n(\mathbf{x})|d\mathbf{x}|.  \label{3.22}
\end{equation}%
Geodesic lines of this metric are orthogonal at any point $\mathbf{x}$ to
the corresponding wave front passing through this point. We assume
everywhere below that the following assumption holds:

\textbf{Assumption 3.1.} \emph{Geodesic lines of metric} (\ref{3.22}) \emph{%
satisfy the regularity condition in }$\mathbb{R}_{+}^{3}$\emph{, i.e. for
each point }$\mathbf{x}\in \mathbb{R}_{+}^{3}$\emph{\ there exists a single
geodesic line }$L(\mathbf{x})$\emph{\ connecting }$\mathbf{x}$\emph{\ with
the plane }$\Sigma $\emph{\ in (\ref{3.200}) and such that }$L(\mathbf{x})$%
\emph{\ intersects the plane }$\Sigma $\emph{\ orthogonally.}

\medskip Below $L(\mathbf{x})$ denotes geodesic lines specified in
Assumption 3.1. The function $\tau (\mathbf{x})$ has the form%
\begin{equation}
\tau (\mathbf{x})=\dint\limits_{L(\mathbf{x})}n(\mathbf{x}^{\prime })d\sigma
,\text{ }\mathbf{x\in }\mathbb{R}_{+}^{3},  \label{3.022}
\end{equation}%
where $\mathbf{x}^{\prime }=\mathbf{x}^{\prime }(\sigma )$ and $d\sigma $ is
the element of the Euclidean length. Denote 
\begin{equation}
\mathbf{p}(\mathbf{x})=\nabla \tau (\mathbf{x}).  \label{3.220}
\end{equation}%
To find $L(\mathbf{x})$ we need to solve the Cauchy problem for the
following system of ordinary differential equations:%
\begin{equation}
\frac{d\mathbf{x}}{ds}=\frac{\mathbf{p}(\mathbf{x})}{n^{2}(\mathbf{x})}%
,\quad \frac{d\mathbf{p}(\mathbf{x})}{ds}=\nabla \ln n(\mathbf{x}),\quad 
\frac{d\tau (\mathbf{x})}{ds}=1,~s>0,  \label{3.23}
\end{equation}%
\begin{equation}
\mathbf{x}|_{s=0}=\mathbf{x}^{0},~\mathbf{p}|_{s=0}=\mathbf{p}^{0},~\tau
|_{s=0}=0,  \label{3.24}
\end{equation}%
where $\mathbf{x}^{0}=(x^{0},y^{0},0)\in \Sigma $. Since the vector $\mathbf{%
p}^{0}$ is orthogonal to the plane $\Sigma $ and $|\mathbf{p}^{0}|=1$ on $%
\Sigma $, then $\mathbf{p}^{0}=(0,0,1).$

For each point $\mathbf{x}^{0}\in \Sigma ,$ the solution of Cauchy problem (%
\ref{3.23}), (\ref{3.24}) defines the geodesic line of the Riemannian metric
(\ref{3.22}), which is orthogonal to the plane $\Sigma $ at the point $%
\mathbf{x}^{0}$. By Assumption 3.1, there exists a one-to-one correspondence
between the points $\mathbf{x\in }\mathbb{R}_{+}^{3}$ and the pairs $\left( 
\mathbf{x}^{0},s\right) \in \Sigma \times \left( 0,\infty \right) $. The
equation of $L(\mathbf{x})$ is given in the form $\mathbf{x}=\mathbf{\xi }(s,%
\mathbf{x}^{0})$. It follows from the last equation (\ref{3.23}) and initial
conditions (\ref{3.24}) that, to find the function $\tau (\mathbf{x}),$ we
need to solve the equation $\mathbf{x}=\mathbf{\xi }(s,\mathbf{x}^{0})$ with
respect to the vector $\left( s,x^{0},y^{0}\right) $. Then we find $s=s(%
\mathbf{x}).$ Next, recalling that by (\ref{3.23}) and (\ref{3.24}) the
parameter $s$ coincides with $\tau (\mathbf{x}),$ we set $\tau (\mathbf{x}%
)=s(\mathbf{x})$. The smoothness of functions $\mathbf{\xi }(s,\mathbf{x}%
^{0})$, $\mathbf{p}(s,\mathbf{x}^{0})$, $s(\mathbf{x})$ is determined by the
smoothness of the function $n(\mathbf{x})$. Since $n\in C^{16}\left( \mathbb{%
R}^{3}\right) $ by (\ref{3.1}), then the function $\tau \in C^{16}\left( 
\mathbb{R}_{+}^{3}\right) $ (see, for instance, [\cite{Rom2}, pp. 26-27],
for the similar derivation).

Introduce the non-negative function $\varphi (\mathbf{x})$ as%
\begin{equation}
\varphi (\mathbf{x})=\left\{ 
\begin{array}{c}
-z\text{ for }z\leq 0, \\ 
\tau (\mathbf{x})\text{ for }z>0.%
\end{array}%
\right.  \label{3.25}
\end{equation}%
Then the function $\varphi \in C^{16}\left( \overline{\mathbb{R}_{+}^{3}}%
\right) .$ The equation $t=\varphi (\mathbf{x})$ defines the characteristic
wedge in $\mathbb{R}^{4}$ for the plane wave originated on the plane $%
\left\{ z=0\right\} $ while it travels inside space $\mathbb{R}^{3}$.

Let $T>0$ be a number. Define the domain $D_{T}$ as 
\begin{equation}
D_{T}=\{(\mathbf{x},t):\,0\leq \varphi (\mathbf{x})\leq t\leq T\}.
\label{3.250}
\end{equation}%
By (\ref{2.2}) the speed $1/n(\mathbf{x})\leq 1.$ Hence, (\ref{3.2100}) and (%
\ref{3.25}) imply that 
\begin{equation}
D_{T}\subset \left\{ (x,y)\in \mathbb{R}^{2},z\in \lbrack -T,T\right\}
\times \left[ 0,T\right] .  \label{3.251}
\end{equation}

\subsection{The structure of the solution of the forward problem}

\label{sec:3.2}

\textbf{Lemma 3.1.} \emph{Let conditions (\ref{2.1})-(\ref{2.3}) be in place
and let} $L(\mathbf{x})$\emph{\ be the geodesic line corresponding to }$\tau
(\mathbf{x})$.\emph{Then the following inequality holds along this geodesic
line:} 
\begin{equation}
\frac{d}{ds}\Delta \tau (\mathbf{x})\leq 6n_{00}^{2},~\mathbf{x}\in \mathbb{R%
}_{+}^{3},  \label{6.02}
\end{equation}%
\emph{where constant $n_{00}$ is defined in (\ref{2.2}).}

\textbf{Proof.} Only in this proof we introduce the following notations: 
\begin{equation*}
\mathbf{x}=(x_{1},x_{2},x_{3}),~\mathbf{p}(\mathbf{x})=\nabla \tau (\mathbf{x%
})\text{, }\tau _{x_{i}x_{j}}(\mathbf{x})=\kappa _{ij}(\mathbf{x}),i,j=1,2,3,
\end{equation*}%
\begin{equation*}
\Delta \tau (\mathbf{x})=\kappa (\mathbf{x})=\sum\limits_{i=1}^{3}\kappa
_{ii}(\mathbf{x}).
\end{equation*}%
Denote $p_{i}=\tau _{x_{i}}$, $i=1,2,3$ and use the eikonal equation (\ref%
{3.21}) $\left\vert \nabla \tau (\mathbf{x})\right\vert ^{2}=n^{2}(\mathbf{x}%
).$ Hence, 
\begin{equation*}
\sum_{i=1}^{3}p_{i}^{2}=n^{2}(\mathbf{x}).
\end{equation*}%
Differentiate this equation with respect to $x_{j}$ and use $%
(p_{i})_{x_{j}}=(p_{j})_{x_{i}}=\kappa _{ij}$. We obtain 
\begin{equation*}
\sum_{i=1}^{3}p_{i}\kappa _{ij}=nn_{x_{j}},~j=1,2,3.
\end{equation*}%
Differentiating this equality with respect to $x_{\ell }$, we get 
\begin{equation*}
\sum_{i=1}^{3}p_{i}(\kappa _{j\ell })_{x_{i}}+\sum_{i=1}^{3}\kappa
_{ij}\kappa _{i\ell }=n_{x_{j}}n_{x_{\ell }}+nn_{x_{\ell }x_{j}},~j,\ell
=1,2,3.
\end{equation*}%
Along the geodesic line $\Gamma (\mathbf{x})$ we can rewrite the latter
equation as 
\begin{equation*}
\frac{d\kappa _{j\ell }}{ds}+\frac{1}{n^{2}}\sum_{i=1}^{3}\kappa _{ij}\kappa
_{i\ell }=\frac{1}{n^{2}}\left( n_{x_{j}}n_{x_{\ell }}+nn_{x_{\ell
}x_{j}}\right) ,~j,\ell =\overline{1,3}.
\end{equation*}%
Take here $\ell =j$ and consider the summation with respect to $j$ from $1$
to $3$. We obtain 
\begin{equation}
\frac{d\kappa }{ds}+\frac{1}{n^{2}(\mathbf{x})}\sum\limits_{i,j=1}^{3}\kappa
_{ij}^{2}=\frac{1}{n^{2}(\mathbf{x})}|\nabla n(\mathbf{x})|^{2}+\frac{1}{n(%
\mathbf{x})}\Delta n(\mathbf{x}).  \label{7.1}
\end{equation}%
It follows from (\ref{2.2}) that 
\begin{equation}
\left\vert \frac{1}{n^{2}(\mathbf{x})}|\nabla n(\mathbf{x})|^{2}+\frac{1}{n(%
\mathbf{x})}\Delta n(\mathbf{x})\right\vert \leq 6n_{00}^{2},~\mathbf{x}\in 
\mathbb{R}_{+}^{3}.  \label{7.100}
\end{equation}%
The target estimate (\ref{6.02}) of Lemma 3.1 follows from (\ref{7.1}) and (%
\ref{7.100}). $\square $

Let $H\left( t\right) $ be the Heaviside function, 
\begin{equation*}
H\left( t\right) =\left\{ 
\begin{array}{c}
1\text{ for }t\geq 0, \\ 
0\text{ for }t<0.%
\end{array}%
\right.
\end{equation*}

\textbf{Theorem 3.1.} \emph{Let conditions (\ref{2.1})-(\ref{2.3}) hold.
Then:}

\emph{1.} \emph{There exists unique solution of problem (\ref{2.12}), (\ref%
{2.13}), which can be represented as } 
\begin{equation}
u(\mathbf{x},t)=H\left( t-\varphi (\mathbf{x})\right) \left[ A(\mathbf{x})+%
\widehat{u}(\mathbf{x},t)\right] ,~\mathbf{x}\in \mathbb{R}^{3},~t\in \left(
0,T\right] ,  \label{3.26}
\end{equation}%
where the function $\widehat{u}\in C^{2}\left( D_{T}\right) $, is compactly
supported in the domain $D_{T}$ and 
\begin{equation}
\lim_{t\rightarrow \varphi (\mathbf{x})^{+}}\widehat{u}(\mathbf{x},t)=0,
\label{3.260}
\end{equation}%
\begin{equation}
A\left( \mathbf{x}\right) >0\text{ in }\mathbb{R}^{3}.  \label{3.261}
\end{equation}%
\emph{The function} $A(\mathbf{x})\in C^{14}\left( \mathbb{R}^{3}\right) $ 
\textit{\ }\emph{and has the form:}%
\begin{equation}
A(\mathbf{x})=\frac{1}{2}\left\{ 
\begin{array}{c}
1\text{ for }z\leq 0, \\ 
\exp \left( -\frac{1}{2}\int\limits_{L(\mathbf{x})}\frac{\Delta \tau (%
\mathbf{x}^{\prime })}{n^{2}(\mathbf{x}^{\prime })}ds\right) \text{ for }z>0,%
\end{array}%
\right.  \label{3.2601}
\end{equation}%
\emph{where} $\mathbf{x}^{\prime }$ \emph{is the variable point along} $L(%
\mathbf{x})$\emph{. }

2. \emph{The inequality }$A\left( \mathbf{x}\right) >0$ \emph{in (\ref{3.261}%
) can be replaced with the following stronger estimate from the below:}%
\begin{equation}
A\left( \mathbf{x}\right) \geq A_{0}=\frac{1}{2}\exp \left(
-3n_{00}^{2}n_{0}^{2}/2\right) ,~\mathbf{x}\in \Omega .  \label{3.262}
\end{equation}%
\emph{\ } \textbf{Proof. } Problem (\ref{2.12}), (\ref{2.13}) is equivalent
with the following one: 
\begin{equation}
\mathcal{L}u\equiv n^{2}\left( \mathbf{x}\right) u_{tt}-\Delta u=0,~\left( 
\mathbf{x},t\right) \in \mathbb{R}^{4},~z\neq 0,~u\mathbf{\mid }_{t<0}=0,
\label{2.120}
\end{equation}%
\begin{equation}
u_{z}|_{z=+0}-u_{z}|_{z=-0}=-\delta \left( t\right)
,~u|_{z=+0}-u|_{z=-0}=0,~(x,y,t)\in \mathbb{R}^{3}.  \label{2.130}
\end{equation}%
Let $r>1$ be an integer, which will be chosen later. We represent the
solution of problem (\ref{2.120}), (\ref{2.130}) in the form 
\begin{equation}
u(\mathbf{x},t)=\sum\limits_{k=0}^{r}\alpha _{k}(\mathbf{x})H_{k}(t-\varphi (%
\mathbf{x}))+u_{r}(\mathbf{x},t),  \label{6.11}
\end{equation}%
where 
\begin{equation}
H_{k}(t)=\frac{t^{k}}{k!}H(t),~k=0,1,2,\ldots ,r.  \label{6.1100}
\end{equation}%
Also, denote $H_{-1}(t)=\delta (t)$ and $H_{-2}(t)=\delta ^{\prime }(t)$.

Recall that the function $\varphi (\mathbf{x})$ is continuous in $\mathbb{R}%
^{3}$, $\varphi (\mathbf{x})=-z$ for $z\leq 0$ and $\varphi \in C^{16}(%
\mathbb{R}_{+}^{3})$. But the derivative $\varphi _{z}(\mathbf{x})$ is
discontinuous across $z=0$, namely, $\varphi _{z}(\mathbf{x})=-1$ for $z\leq
0$ and $\varphi _{z}(\mathbf{x})=1$ for $z=+0$. Thus, 
\begin{equation}
\varphi |_{z=0}=0,~\varphi _{z}|_{z=+0}=1,~\varphi _{z}|_{z=-0}=-1.
\label{6.1200}
\end{equation}%
Taking into account (\ref{6.1200}), we need to consider representation (\ref%
{6.11}) separately for $z<0$ and $z\geq 0$.

We seek functions $\alpha _{k}(\mathbf{x})$ in the form: 
\begin{equation}
\alpha _{k}(\mathbf{x})=\left\{ 
\begin{array}{ll}
\alpha _{k}^{-}(\mathbf{x}), & z<0, \\ 
\alpha _{k}^{+}(\mathbf{x}), & z\geq 0.%
\end{array}%
\right.  \label{6.12}
\end{equation}%
Substituting representation (\ref{6.11}) in (\ref{2.130}) and equating
coefficients at $H_{k}(t)$, we obtain 
\begin{equation}
\left\{ 
\begin{array}{l}
-(\alpha _{k}\varphi _{z})_{z=+0}+(\alpha _{k}\varphi _{z})_{z=-0}=-\delta
_{k0}, \\ 
\alpha _{k}|_{z=+0}-\alpha _{k}|_{z=-0}=0,%
\end{array}%
\right.  \label{6.120}
\end{equation}%
where $\delta _{k0}$ is the Kronecker's delta. Then function $u_{r}(\mathbf{x%
},t)$ satisfies the following conjugate conditions 
\begin{equation}
(u_{r})_{z}|_{z=+0}-(u_{r})_{z}|_{z=-0}=0,~u_{r}|_{z=+0}-u_{r}|_{z=-0}=0.
\label{6.121}
\end{equation}%
Using (\ref{6.1200}), we find from equations (\ref{6.120}) 
\begin{equation}
\alpha _{k}^{-}|_{z=-0}=\alpha _{k}^{+}|_{z=0}=\frac{\delta _{k0}}{2}%
,~k=0,1,\ldots ,r.  \label{6.122}
\end{equation}%
Apply the operator $\mathcal{L}u$ for $z\neq 0$ to both sides of (\ref{6.11}%
).\ First, 
\begin{equation}
\partial _{t}^{2}\left[ n^{2}\left( \mathbf{x}\right) \alpha _{k}(\mathbf{x}%
)H_{k}(t-\varphi (\mathbf{x}))\right] =\alpha _{k}(\mathbf{x})n^{2}\left( 
\mathbf{x}\right) H_{k-2}(t-\varphi (\mathbf{x}).  \label{40}
\end{equation}%
Second,%
\begin{equation*}
-\Delta \left[ \alpha _{k}(\mathbf{x})H_{k}(t-\varphi (\mathbf{x}))\right]
=-H_{k}(t-\varphi (\mathbf{x}))\Delta \alpha _{k}(\mathbf{x})-2\nabla \alpha
_{k}(\mathbf{x})\cdot \nabla H_{k}(t-\varphi (\mathbf{x}))
\end{equation*}%
\begin{equation*}
-\alpha _{k}(\mathbf{x})\Delta \left[ H_{k}(t-\varphi (\mathbf{x}))\right]
\end{equation*}%
\begin{equation}
=-H_{k}(t-\varphi (\mathbf{x}))\Delta \alpha _{k}(\mathbf{x})+2\nabla \alpha
_{k}(\mathbf{x})\cdot \nabla \varphi (\mathbf{x})H_{k-1}(t-\varphi (\mathbf{x%
}))  \label{41}
\end{equation}%
\begin{equation*}
+\alpha _{k}(\mathbf{x})\Delta \varphi \left( \mathbf{x}\right)
H_{k-1}(t-\varphi (\mathbf{x}))-\alpha _{k}(\mathbf{x})\left\vert \nabla
\varphi \left( \mathbf{x}\right) \right\vert ^{2}H_{k-2}(t-\varphi (\mathbf{x%
})).
\end{equation*}%
Since by (\ref{3.21}) $n^{2}\left( \mathbf{x}\right) -\left\vert \nabla
\varphi \left( \mathbf{x}\right) \right\vert ^{2}=0,$ then (\ref{40}) and (%
\ref{41}) imply: 
\begin{equation*}
\left( n^{2}\left( \mathbf{x}\right) \partial _{t}^{2}-\Delta \right) \left[
\alpha _{k}(\mathbf{x})H_{k}(t-\varphi (\mathbf{x}))\right] =
\end{equation*}%
\begin{equation*}
-H_{k}(t-\varphi (\mathbf{x}))\Delta \alpha _{k}(\mathbf{x})+2\nabla \alpha
_{k}(\mathbf{x})\nabla \varphi (\mathbf{x})H_{k-1}(t-\varphi (\mathbf{x}%
))+\alpha _{k}(\mathbf{x})\Delta \varphi \left( \mathbf{x}\right)
H_{k-1}(t-\varphi (\mathbf{x})).
\end{equation*}%
Hence, 
\begin{equation*}
\sum\limits_{k=0}^{r}\alpha _{k}(\mathbf{x})H_{k}(t-\varphi (\mathbf{x}%
))=\sum\limits_{k=0}^{r}\left[ 2\nabla \alpha _{k}(\mathbf{x})\nabla \varphi
(\mathbf{x})+\alpha _{k}(\mathbf{x})\Delta \varphi \left( \mathbf{x}\right) %
\right] H_{k-1}(t-\varphi (\mathbf{x}))
\end{equation*}%
\begin{equation*}
-\sum\limits_{k=0}^{r}[\Delta \alpha _{k}(\mathbf{x})]H_{k}(t-\varphi (%
\mathbf{x})).
\end{equation*}%
Next, 
\begin{equation*}
-\sum\limits_{k=0}^{r}[\Delta \alpha _{k}(\mathbf{x})]H_{k}(t-\varphi (%
\mathbf{x}))=-\sum\limits_{k=1}^{r+1}[\Delta \alpha _{k-1}(\mathbf{x}%
)]H_{k-1}(t-\varphi (\mathbf{x}))=
\end{equation*}%
\begin{equation*}
=-\sum\limits_{k=0}^{r}[\Delta \alpha _{k-1}(\mathbf{x})]H_{k-1}(t-\varphi (%
\mathbf{x}))-[\Delta \alpha _{r}(\mathbf{x})]H_{r}(t-\varphi (\mathbf{x})),
\end{equation*}%
where we formally set $\alpha _{-1}(\mathbf{x})\equiv 0$. Hence, we obtain%
\begin{equation*}
\mathcal{L}u=\sum\limits_{k=0}^{r}[2\nabla \alpha _{k}(\mathbf{x})\cdot
\nabla \varphi (\mathbf{x})+\alpha _{k}(\mathbf{x})\Delta \varphi (\mathbf{x}%
)-\Delta \alpha _{k-1}(\mathbf{x})]H_{k-1}(t-\varphi (\mathbf{x}))
\end{equation*}%
\begin{equation*}
-[\Delta \alpha _{r}(\mathbf{x})]H_{r}(t-\varphi (\mathbf{x}))+\mathcal{L}%
u_{r}(\mathbf{x},t).~
\end{equation*}

Equating here to zero terms at $H_{k-1}(t-\varphi (\mathbf{x}))$ and taking
into account conditions (\ref{6.122}), we obtain equations for $\alpha
_{k}^{-}$ and $\alpha _{k}^{+}$, $k=0,1,\ldots ,r$, 
\begin{equation}
2\nabla \alpha _{k}^{-}(\mathbf{x})\cdot \nabla \varphi(\mathbf{x})+\alpha
_{k}^{-}(\mathbf{x})\Delta \varphi(\mathbf{x})=\Delta \alpha _{k-1}^{-}(%
\mathbf{x}),~z<0, ~\alpha _{k}^{-}|_{z=-0}=\frac{1}{2}\delta _{k0},
\label{6.131}
\end{equation}%
\begin{equation}
2\nabla \alpha _{k}^{+}(\mathbf{x})\cdot \nabla \tau(\mathbf{x})+\alpha
_{k}^{+}(\mathbf{x})\Delta \tau(\mathbf{x})=\Delta \alpha _{k-1}^{+}(\mathbf{%
x}),~z>0,~\alpha _{k}^{+}|_{z=0}=\frac{1}{2}\delta _{k0},  \label{6.13}
\end{equation}%
where $\alpha _{-1}^{-}(\mathbf{x})=\alpha _{-1}^{+}(\mathbf{x})=0$.

Note that in equations (\ref{6.131}) $\varphi (\mathbf{x})=-z$. Therefore it
can be written as follows 
\begin{equation}
-2(\alpha _{k}^{-})_{z}(\mathbf{x})=\Delta \alpha _{k-1}^{-}(\mathbf{x}%
),~z<0,~\alpha _{k}^{-}|_{z=0}=\frac{1}{2}\delta _{k0},~k=0,1,\ldots ,r.
\label{6.1310}
\end{equation}%
It follows from (\ref{6.1310}) that 
\begin{equation}
\alpha _{k}^{-}(\mathbf{x})=\frac{1}{2}\delta _{k0},~k=0,1,\ldots ,r.
\label{6.13101}
\end{equation}%
Since $\mathcal{L}u(\mathbf{x},t)=0$ for $z\neq 0$ and conditions (\ref%
{6.121}) hold, the equation for the function $u_{r}(\mathbf{x},t)$ is: 
\begin{equation}
n^{2}(\mathbf{x})(u_{r})_{tt}-\Delta u_{r}=F_{r}(\mathbf{x},t),~(\mathbf{x}%
,t)\in \mathbb{R}^{4},\quad u_{r}|_{t<0}=0,  \label{6.14}
\end{equation}%
where 
\begin{equation}
F_{r}(\mathbf{x},t)=H_{r}(t-\varphi (\mathbf{x}))\Delta \alpha _{r}(\mathbf{x%
}).  \label{6.140}
\end{equation}

Integrate now equation (\ref{6.13}) along the geodesic line $L(\mathbf{x})$.
By (\ref{3.220}) and (\ref{3.23}) we have along this line $\nabla \tau(%
\mathbf{x})=\mathbf{p}(\mathbf{x})=n^{2}(\mathbf{x})d\mathbf{x}/ds$.
Therefore, 
\begin{equation*}
\nabla \alpha _{k}^{+}(\mathbf{x})\cdot \nabla \tau(\mathbf{x})=n^{2}(%
\mathbf{x})\nabla \alpha _{k}^{+}(\mathbf{x})\cdot \frac{d\mathbf{x}}{ds}%
=n^{2}(\mathbf{x})\frac{d}{ds}\alpha _{k}^{+}(\mathbf{x}).
\end{equation*}%
Hence, (\ref{6.13}) is equivalent with: 
\begin{equation}
2n^{2}(\mathbf{x})\frac{d}{ds}\alpha _{k}^{+}(\mathbf{x})+\alpha _{k}^{+}(%
\mathbf{x})\Delta \tau(\mathbf{x})=\Delta \alpha _{k-1}^{+}(\mathbf{x}%
),~\alpha _{k}^{+}|_{z=0}=\frac{1}{2}\delta _{0k}.  \label{6.15}
\end{equation}%
The solution of the Cauchy problem 
\begin{equation*}
2n^{2}(\mathbf{x})\frac{d}{ds}\alpha _{0}^{+}(\mathbf{x})+\alpha _{0}^{+}(%
\mathbf{x})\Delta \tau(\mathbf{x})=0,~\alpha _{0}^{+}|_{z=0}=\frac{1}{2}
\end{equation*}%
is given by the formula 
\begin{equation}
\alpha _{0}^{+}(\mathbf{x})=\frac{1}{2}\exp \left( -\frac{1}{2}%
\int\limits_{L(\mathbf{x})}\frac{\Delta \tau(\mathbf{x}^{\prime })}{n^{2}(%
\mathbf{x}^{\prime })}ds\right) ,  \label{6.16}
\end{equation}%
where $\mathbf{x}^{\prime }$ is a variable point along $L(\mathbf{x})$.
Next, dividing both sides of equation (\ref{6.15}) by $2n^{2}(\mathbf{x})$
we can rewrite (\ref{6.15}) in the form 
\begin{equation}
\exp \left( -\frac{1}{2}\int\limits_{L(\mathbf{x})}\frac{\Delta \tau(\mathbf{%
x}^{\prime })}{n^{2}(\mathbf{x}^{\prime })}ds\right) \frac{d}{ds}\left[
\alpha _{k}^{+}(\mathbf{x})\exp \left( \frac{1}{2}\int\limits_{L(\mathbf{x})}%
\frac{\Delta \tau(\mathbf{x}^{\prime })}{n^{2}(\mathbf{x}^{\prime })}%
ds\right) \right] =  \label{6.1600}
\end{equation}%
\begin{equation*}
=\frac{\Delta \alpha _{k-1}^{+}(\mathbf{x})}{2n^{2}(\mathbf{x})}.
\end{equation*}%
It follows from (\ref{6.15}) and (\ref{6.16}) that (\ref{6.1600}) is
equivalent with 
\begin{equation}
\frac{d}{ds}\left( \frac{\alpha _{k}^{+}(\mathbf{x})}{\alpha _{0}^{+}(%
\mathbf{x})}\right) =\frac{\Delta \alpha _{k-1}^{+}(\mathbf{x})}{2n^{2}(%
\mathbf{x})\alpha _{0}^{+}(\mathbf{x})},~\alpha
_{k}^{+}|_{z=0}=0,~k=1,\ldots ,r.  \label{6.17}
\end{equation}%
Integrating (\ref{6.17}) with respect to $s$ along $L(\mathbf{x})$, we
obtain 
\begin{equation}
\alpha _{k}^{+}(\mathbf{x})=\frac{1}{2}\alpha _{0}^{+}(\mathbf{x}%
)\int\limits_{L(\mathbf{x})}\frac{\Delta \alpha _{k-1}^{+}(\mathbf{x}%
^{\prime })}{n^{2}(\mathbf{x}^{\prime })\alpha _{0}^{+}(\mathbf{x}^{\prime })%
}ds,~k=1,\ldots ,r.  \label{6.18}
\end{equation}

Let $m$ be a sufficiently large integer, which we will choose below. If the
function $n\in C^{m}(\mathbb{R}^{3})$, then 
\begin{equation*}
\tau \in C^{m}(\mathbb{R}_{+}^{3}),~\alpha _{k}^{+}\in C^{m-2k-2}(\mathbb{R}%
_{+}^{3}),~\Delta \alpha _{r}^{+}\in C^{m-2r-4}(\mathbb{R}_{+}^{3}).
\end{equation*}%
Define the domain $G_{X}$ as%
\begin{equation}
G_{X}=\{\mathbf{x}\in \mathbb{R}_{+}^{3}:\,|x|\geq X,|y|\geq X,z\geq 0\}.
\label{6.180}
\end{equation}%
It follows from (\ref{6.16}) and (\ref{6.18}) that $\alpha _{0}^{+}(\mathbf{x%
})=1/2$ and functions $\alpha _{k}^{+}(\mathbf{x})=0$, $k=1,\ldots ,r$, for $%
\mathbf{x}\in G_{X}$. Indeed, since by (\ref{2.3}) $n\left( \mathbf{x}%
\right) =1$ in $G_{X},$ then Assumption 3.1 implies that if $\mathbf{x}\in
G_{X},$ then the geodesic line $L(\mathbf{x})$ is a segment of the straight
line orthogonal to the plane $z=0,$ and this line does not intersect $\Omega 
$.

Setting $A(\mathbf{x})=\alpha _{0}(\mathbf{x})$ and using (\ref{6.11}), (\ref%
{6.12}), (\ref{6.13101}) and (\ref{6.16}), we obtain (\ref{3.2601}).

It follows from (\ref{6.140}) and the above arguments that the function $%
F_{r}(\mathbf{x},t)$ possesses the following properties:

1. For any $T>0$ 
\begin{equation*}
support\left( F_{r}(\mathbf{x},t)\right) \subset D_{T}^{+}=\{(\mathbf{x}%
,t)\in D_{T}:\,z\geq 0\},
\end{equation*}
where the set $D_{T}$ is defined in (\ref{3.250}).

2. $F_{r}(\mathbf{x},t)=0$ for $\mathbf{x}\in G_{X}$.

Furthermore, the projection of the set $D_{T}^{+}$ on the space $\mathbb{R}%
^{3}$ coincides with the set $Y_{T}=\{\mathbf{x}\in \mathbb{R}^{3}:\,0\leq
\tau (\mathbf{x})\leq T\}$. It follows from (\ref{3.251}) that this set is
bounded with respect to $z$, i.e. $Y_{T}\subset \left\{ (x,y)\in \mathbb{R}%
^{2},z\in \lbrack 0,T]\right\} $. Hence, the set $\left( Y_{T}\setminus
G_{X}\right) \subset \{\left\vert x\right\vert ,\left\vert y\right\vert
<X,0\leq z\leq T\}$ is bounded, where $G_{X}$ is the set defined in (\ref%
{6.180}). Thus, the function $F_{r}(\mathbf{x},t)=0$ outside of the finite
domain\ $\{\left\vert x\right\vert ,\left\vert y\right\vert <X,0\leq z\leq
T\}$, i.e. it is compactly supported in $\mathbb{R}_{T}^{4}$. Moreover, $%
F_{r}\in H^{d}(\mathbb{R}_{T}^{4})$, where $d=\min (m-2r-4,r)$. Using the
general theory of hyperbolic equations \cite[Chapter 4]{Lad}, we conclude
that the unique solution $u_{r}\in H^{d+1}(\mathbb{R}_{T}^{4})$ of the
Cauchy problem (\ref{6.14}), (\ref{6.140}) exists and this solution is also
compactly supported in $\mathbb{R}_{T}^{4}$. Moreover, since $F_{r}(\mathbf{x%
},t)=0$ for $(\mathbf{x},t)\notin D_{T}$, then $u_{r}(\mathbf{x},t)=0$ for $%
t<\varphi (\mathbf{x}).$ Embedding theorem implies $u_{r}\in C^{2}\left( 
\overline{\mathbb{R}_{T}^{4}}\right) $ if $d+1>4$. Choose $d=r=4$ and $m=16$%
. Then $u_{r}\in C^{2}\left( \overline{\mathbb{R}_{T}^{4}}\right) $. Since%
\textbf{\ }$u_{r}(\mathbf{x},t)=0$\textbf{\ }for\textbf{\ }$t$\textbf{$<$}$%
\varphi (\mathbf{x})$\textbf{, }we conclude that\textbf{\ }$u_{r}\in
C^{2}\left( \overline{D_{T}}\right) $ and $u_{r}(\mathbf{x},t)=0$ for $%
t=\varphi (\mathbf{x})$ together with derivatives up to the second order.
This explains our smoothness condition (\ref{3.1}).

Finally, setting in (\ref{6.11}) 
\begin{equation*}
A(\mathbf{x})=\alpha_0(\mathbf{x}), ~ \widehat{u}(\mathbf{x}%
,t)=\sum\limits_{k=1}^{r}\alpha _{k}(\mathbf{x})\frac{(t-\varphi (\mathbf{x}%
))^{k}}{k!}+u_{r}(\mathbf{x},t),
\end{equation*}%
we finish the proof of (\ref{3.26}) and (\ref{3.260}).

We now want to prove that $A(\mathbf{x})=\alpha _{0}(\mathbf{x})\geq A_{0}$
with the positive number $A_{0}$ defined in (\ref{3.262}). We use Lemma 3.1
for this purpose.

The function $\alpha _{0}(\mathbf{x})=1/2$ for $z\leq 0$. Using (\ref{6.02})
and (\ref{6.16}), estimate now the function $\alpha _{0}^{+}(\mathbf{x})$.
Note first that 
\begin{equation}
\Delta \tau (\mathbf{x})|_{z=0}=0.  \label{6.160}
\end{equation}%
Indeed, eikonal equation (\ref{3.21}) implies that 
\begin{equation}
\nabla \partial _{z}\tau (\mathbf{x})\cdot \nabla \tau (\mathbf{x})=n(%
\mathbf{x})n_{z}(\mathbf{x}).  \label{100}
\end{equation}%
It follows from the condition $\tau (\mathbf{x})|_{z=0}=0$ in (\ref{3.21})
that, at $z=0$, $(\tau )_{xx}=(\tau )_{yy}=0$, $(\tau )_{z}=1$ and $(\tau
)_{xz}=(\tau )_{yz}=0$. Therefore (\ref{100}) at $z=0$ becomes 
\begin{equation*}
(\tau )_{zz}(x,y,0)=n(x,y,0)n_{z}(x,y,0)=0.
\end{equation*}%
Here the equality $n_{z}(x,y,0)=0$ follows from (\ref{3.1}) and (\ref{2.3}).
Thus, (\ref{6.160}) holds. Hence, integrating inequality (\ref{6.02}) of
Lemma 3.1 with respect to $s\in \left( 0,s^{\prime }\right) $, we conclude
that along $L(\mathbf{x})$ 
\begin{equation}
\frac{\Delta \tau (\mathbf{x}^{\prime })}{2n^{2}(\mathbf{x}^{\prime })}\leq
3n_{00}^{2}s^{\prime },~\mathbf{x}^{\prime }\in \mathbb{R}_{+}^{3},
\label{6.19}
\end{equation}%
Formulae (\ref{6.16}) and (\ref{6.19}) imply: 
\begin{equation*}
\alpha _{0}^{+}(\mathbf{x})\geq \frac{1}{2}\exp \left( -\int\limits_{\Gamma (%
\mathbf{x})}3n_{00}^{2}s^{\prime }ds^{\prime }\right) =\frac{1}{2}\exp
\left( -\frac{3n_{00}^{2}}{2}\tau ^{2}(\mathbf{x})\right) ,~\mathbf{x}\in 
\mathbb{R}_{+}^{3},
\end{equation*}%
It follows from (\ref{2.2}) and eikonal equation (\ref{3.21}) that 
\begin{equation}
\partial _{z}\tau \leq n_{0}\text{ \ in }\mathbb{R}_{+}^{3}.  \label{6.190}
\end{equation}
Since $\tau (x,y,0)=0,$ then 
\begin{equation*}
\tau (\mathbf{x})=\int\limits_{0}^{z}\partial _{r}\tau (x,y,r)dr\leq n_{0},~%
\mathbf{x}\in \overline{\Omega }.
\end{equation*}%
Hence, 
\begin{equation*}
A(\mathbf{x})\geq \alpha _{0}^{+}(\mathbf{x})\geq \ \frac{1}{2}\exp \left( -%
\frac{3n_{00}^{2}}{2}n_{0}^{2}\right) =A_{0},~\mathbf{x}\in \overline{\Omega 
}.\text{ }~~\square
\end{equation*}

\section{A Boundary Value Problem in Partial Finite Differences}

\label{sec:4}

\subsection{The boundary value problem}

\label{sec:4.1}

Consider the function $v\left( \mathbf{x},t\right) ,$ 
\begin{equation}
v\left( \mathbf{x},t\right) =\dint\limits_{0}^{t}u(\mathbf{x},s)ds.
\label{3.35}
\end{equation}%
Recall that by (\ref{3.25}) $\varphi (\mathbf{x})=\tau (\mathbf{x})$ for $%
\mathbf{x}\in \Omega .$ Hence, by (\ref{3.26}), (\ref{3.260}) and (\ref{3.35}%
) 
\begin{equation}
v(\mathbf{x},t)=\left( t-\tau (\mathbf{x})\right) H\left( t-\tau (\mathbf{x}%
)\right) \left[ A(\mathbf{x})+\widehat{v}(\mathbf{x},t)\right] ,~\mathbf{x}%
\in \Omega ,  \label{3.36}
\end{equation}%
where $\widehat{v}(\mathbf{x},\tau (\mathbf{x})+0)=0$.

Estimate $\max_{\mathbf{x}\in \overline{\Omega }}\tau (\mathbf{x})$. Since $%
\tau |_{z=0}=0,$ then, using (\ref{6.190}), we obtain 
\begin{equation*}
\tau \left( x,y,z\right) =\dint\limits_{0}^{z}\partial _{r}\tau \left(
x,y,r\right) dr\leq n_{0}z\leq n_{0},\text{ }\mathbf{x=}\left( x,y,z\right)
\in \Omega .
\end{equation*}%
Hence, 
\begin{equation}
\max_{\mathbf{x}\in \overline{\Omega }}\tau (\mathbf{x})\leq n_{0}
\label{5.970}
\end{equation}%
We assume that $T>n_{0}$. Denote 
\begin{equation}
T_{1}=T-n_{0}>0.  \label{3.34}
\end{equation}%
Also, denote $Q_{T_{1}}=\Omega \times \left( 0,T_{1}\right) .$ Let 
\begin{equation}
P\left( \mathbf{x},t\right) =v(\mathbf{x},t+\tau (\mathbf{x}))\text{ for }(%
\mathbf{x},t)\in Q_{T_{1}}.  \label{3.37}
\end{equation}%
This function is defined for all $(\mathbf{x},t)\in Q_{T_{1}}$ since by (\ref%
{5.970}) 
\begin{equation*}
0<t+\tau (\mathbf{x})<T_{1}+\max_{\mathbf{x}\in \overline{\Omega }}\tau (%
\mathbf{x})\leq T_{1}+n_{0}=T\text{ \ in }Q_{T_{1}}.
\end{equation*}%
Hence, by (\ref{3.36}), (\ref{3.34}) and (\ref{3.37})%
\begin{equation}
P\left( \mathbf{x},t\right) =t\left[ A(\mathbf{x})+\widehat{P}(\mathbf{x},t)%
\right] ,\text{ }\left( \mathbf{x},t\right) \in Q_{T_{1}},  \label{3.38}
\end{equation}%
\begin{equation}
\lim_{t\rightarrow 0^{+}}\widehat{P}(\mathbf{x},t)=0.  \label{3.380}
\end{equation}%
Using (\ref{2.12}), (\ref{3.21}), (\ref{3.35}) and (\ref{3.37}), we obtain 
\begin{equation}
\Delta P-2\nabla_{\mathbf{x}}P_{t}\cdot\nabla\tau -P_{t}\Delta \tau =0,\text{
}\left( \mathbf{x,}t\right) \in Q_{T_{1}}.  \label{3.39}
\end{equation}%
Denote 
\begin{equation}
w(\mathbf{x},t)=P_{t}\left( \mathbf{x},t\right) .  \label{3.390}
\end{equation}%
Note that by (\ref{3.262}), (\ref{3.38}) and (\ref{3.380}) 
\begin{equation}
w(\mathbf{x},0)=A\left( \mathbf{x}\right) \geq A_{0}=\frac{1}{2}\exp \left(
-3n_{00}^{2}n_{0}^{2}/2\right) >0,~\mathbf{x}\in \Omega .  \label{3.41}
\end{equation}%
Setting in (\ref{3.39}) $t=0$, using the fact that by (\ref{3.38}) and (\ref%
{3.380}) $\Delta P\left( x,0\right) =0$ and also using (\ref{3.390}) and (%
\ref{3.41}), we obtain $\Delta \tau +2\left( \nabla \ln w\left( \mathbf{x}%
,0\right) \right) \cdot \nabla \tau =0$ for $\mathbf{x}\in \Omega .$ Hence,
using this equation and (\ref{3.38}), we obtain the $2\times 2$ system of
non local nonlinear PDEs: 
\begin{equation}
\left\{ 
\begin{array}{c}
\Delta w-2\nabla_{\mathbf{x}}w_{t}\cdot\nabla\tau +2w_{t}\left( \nabla \ln
w\left( \mathbf{x},0\right) \right) \cdot \nabla \tau =0,\text{ }\left( 
\mathbf{x,}t\right) \in Q_{T_{1}}, \\ 
\Delta \tau +2\left( \nabla \ln w\left( \mathbf{x},0\right) \right) \cdot
\nabla \tau =0,\text{ }\mathbf{x}\in \Omega .%
\end{array}%
\right.  \label{3.40}
\end{equation}

To find boundary conditions for system (\ref{3.40}), we use functions $%
f_{0}\left( \boldsymbol{x},t\right) ,f_{1}\left( \mathbf{x},t\right)
,f_{2}\left( \mathbf{x},t\right) $ in (\ref{2.14}). Let%
\begin{equation*}
g_{0}\left( \mathbf{x},t\right) =f_{0}\left( \mathbf{x},\tau \left( \mathbf{x%
}\right) +t\right) ,\text{ }\left( \mathbf{x},t\right) \in \Gamma _{T_{1}},
\end{equation*}%
\begin{equation*}
g_{1}\left( \mathbf{x},t\right) =f_{1}\left( \mathbf{x},\tau \left( \mathbf{x%
}\right) +t\right) +\partial _{t}g_{0}\left( \mathbf{x},t\right) \partial
_{z}\tau \left( \mathbf{x}\right) ,~\left( \mathbf{x},t\right) \in \Gamma
_{T_{1}},
\end{equation*}%
\begin{equation*}
g_{2}\left( \mathbf{x},t\right) =f_{2}\left( \mathbf{x},\tau \left( \mathbf{x%
}\right) +t\right) ,~\left( \mathbf{x},t\right) \in \Theta _{T_{1}}.
\end{equation*}

Since by (\ref{2.6}) $\Gamma \subset \left\{ z=0\right\} $ and since by (\ref%
{2.3}) $n\left( \mathbf{x}\right) \mid _{\Gamma }=1$ and also since by the
second condition in (\ref{3.21}) $\tau \mid _{\Gamma }=0,$ then $\partial
_{x}\tau \mid _{\Gamma }=\partial _{y}\tau \mid _{\Gamma }=0.$ Hence, (\ref%
{3.21}) and (\ref{3.210}) imply $\partial _{z}\tau \mid _{\Gamma }=n\left( 
\mathbf{x}\right) \mid _{\Gamma }=1.$ Also, it obviously follows from (\ref%
{2.70}) and (\ref{2.3}) that $\tau \left( \mathbf{x}\right) \mid _{\Theta
}=z.$ Thus, the boundary conditions for system (\ref{3.40}) are:%
\begin{equation}
w\mid _{\Gamma _{T_{1}}}=g_{0}\left( \mathbf{x},t\right) ,\text{ }w_{z}\mid
_{\Gamma _{T_{1}}}=g_{1}\left( \mathbf{x},t\right) ,w\mid _{\Theta
_{T_{1}}}=g_{2}\left( \mathbf{x},t\right) ,  \label{3.45}
\end{equation}%
\begin{equation}
\tau \mid _{\Gamma }=0,~\partial _{z}\tau \mid _{\Gamma }=1,~\tau \mid
_{\Theta }=z.  \label{3.46}
\end{equation}

Therefore, we arrive at the following Boundary Value Problem for the $%
2\times 2$ system (\ref{3.40}) of nonlinear and non local PDEs:

\textbf{Boundary Value Problem } (BVP): \emph{Find the pair of functions }$%
\left( w,\tau \right) \in C^{2}\left( \overline{Q}_{T_{1}}\right) \times
C^{2}\left( \overline{\Omega }\right) $\emph{\ satisfying equations (\ref%
{3.40}) and boundary conditions (\ref{3.45}), (\ref{3.46}), assuming that (%
\ref{3.41}) holds. }

As it was pointed out in Introduction, we cannot prove stability estimates
for this BVP. However, we can prove the desired stability estimates if we
rewrite this BVP in the form of partial finite differences, in which the
derivatives with respect to $x$ and $y$ are written in finite differences,
whereas the derivatives with $z$ and $t$ are written in the conventional
continuous way. In doing so, we assume that the step size $h$ of the finite
difference scheme is bounded from the below by a fixed positive number. The
latter assumption is a quite natural one in computations. Thus, we rewrite
in subsection 4.2 the above BVP in partial finite differences.

\subsection{Partial finite differences}

\label{sec:4.2}

For brevity, we use the same grid step size $h$ in both $x$ and $y$
directions. Choose a small number $h_{0}\in \left( 0,1\right) .$ We assume
everywhere below in this paper that 
\begin{equation}
h\geq h_{0}.  \label{4.1}
\end{equation}%
Consider two partitions of the interval $\left[ -X,X\right] ,$ 
\begin{align*}
& -X=x_{0}<x_{1}<\ldots <x_{N}<x_{N+1}=X,\quad x_{i}-x_{i-1}=h, \\
& -X=y_{0}<y_{1}<\ldots <y_{N}<y_{N+1}=X,\quad y_{i}-y_{i-1}=h.
\end{align*}%
Thus, the interior grid points are $\left\{ \left( x_{i},y_{j}\right)
\right\} _{i,j=1}^{N}$ and other grid points are located on the part of the
boundary $\Theta \subset \partial \Omega ,$ see (\ref{2.70}). Let $\mathbf{x}%
^{h}=\left( x_{i},y_{j},z\right) _{i,j=0}^{N+1}$ denotes semi-discrete
points. For any function $s\left( \mathbf{x},t\right) $ defined on the set $%
\overline{Q_{T_{1}}},$ we denote by $s\left( \mathbf{x}^{h},t\right) $ the
corresponding semi-discrete vector function defined at points $\left( 
\mathbf{x}^{h},t\right) \in \overline{\Omega }\times \left[ 0,T_{1}\right] .$
We also introduce the following notations: 
\begin{equation*}
\Omega _{h}=\left\{ \mathbf{x}^{h}=\left( x_{i},y_{j},z\right)
_{i,j=1}^{N},z\in \left( 0,1\right) \right\} ,
\end{equation*}%
\begin{equation*}
\overline{\Omega }_{h}=\left\{ \mathbf{x}^{h}=\left( x_{i},y_{j},z\right)
_{i,j=0}^{N+1},z\in \left[ 0,1\right] \right\} ,
\end{equation*}%
\begin{equation*}
Q_{h,T_{1}}=\Omega _{h}\times \left( 0,T_{1}\right) ,\text{ }\overline{%
Q_{h,T_{1}}}=\overline{\Omega }_{h}\times \left[ 0,T_{1}\right] ,
\end{equation*}%
\begin{equation*}
\Gamma _{h}=\left\{ \mathbf{x}^{h}=\left\{ \left( x_{i},y_{j},0\right)
\right\} _{i,j=1}^{N}\right\} ,\text{ }\Gamma _{h,T_{1}}=\Gamma _{h}\times
\left( 0,T_{1}\right) ,
\end{equation*}%
\begin{equation*}
\Gamma _{h}^{\prime }=\left\{ \mathbf{x}^{h}=\left\{ \left(
x_{i},y_{j},1\right) \right\} _{i,j=1}^{N}\right\} ,\text{ }\Gamma
_{h,T_{1}}^{\prime }=\Gamma _{h}^{\prime }\times \left( 0,T_{1}\right) ,
\end{equation*}%
\begin{equation*}
\Theta _{h}=\left\{ \mathbf{x}^{h}\in \Theta \right\} ,\text{ }\Theta
_{h,T_{1}}=\Theta _{h}\times \left( 0,T_{1}\right) .
\end{equation*}%
The derivatives in finite differences with respect to $x$ are defined as:%
\begin{equation}
\partial _{x}s_{i,j}^{h}\left( z\right) =\frac{s_{i-1,j}^{h}\left(
z,t\right) -s_{i+1,j}^{h}\left( z,t\right) }{2h},\text{ }i,j=1,...,N,
\label{4.2}
\end{equation}%
\begin{equation}
\partial _{x}s^{h}\left( \mathbf{x}^{h},t\right) =s_{x}^{h}\left( \mathbf{x}%
^{h},t\right) =\left\{ \partial _{x}s_{i,j}^{h}\left( z,t\right) \right\}
_{i,j=1}^{N},  \label{4.4}
\end{equation}%
\begin{equation}
\partial _{x}^{2}s_{i,j}^{h}\left( z,t\right) =\frac{s_{i+1,j}^{h}\left(
z,t\right) -2s_{i,j}^{h}\left( z,t\right) +s_{i-1,j}^{h}\left( z,t\right) }{%
h^{2}},\text{ }i,j=1,...,N,  \label{4.5}
\end{equation}%
\begin{equation}
\partial _{x}^{2}s^{h}\left( \mathbf{x}^{h},t\right) =s_{xx}^{h}\left( 
\mathbf{x}^{h},t\right) =\left\{ \partial _{x}^{2}s_{i,j}^{h}\left(
z,t\right) \right\} _{i,j=1}^{N}.  \label{4.8}
\end{equation}

Derivatives $s_{y}^{h},s_{yy}^{h}$ are defined completely similarly.
Formulas (\ref{4.2})-(\ref{4.8}) as well as their analogs for $%
s_{y}^{h},s_{yy}^{h}$ fully define finite difference versions of $x,y-$%
derivatives involved in equations (\ref{3.40}). Everywhere below the
corresponding Laplace operator in the partial finite differences as well as
the gradient vector are given by: 
\begin{equation}
\Delta ^{h}s^{h}=s_{xx}^{h}+s_{yy}^{h}+s_{zz}^{h},~\nabla ^{h}s^{h}=\left(
s_{x}^{h},s_{y}^{h},s_{z}^{h}\right) ,  \label{4.80}
\end{equation}%
where the $z-$derivatives are understood in the regular manner, and the same
for the $t-$derivatives in follow up formulas.

We introduce the semi-discrete analogs of conventional function spaces of
real valued functions as:%
\begin{equation*}
H^{2,h}\left( Q_{h,T_{1}}\right) =\left\{ 
\begin{array}{c}
s^{h}\left( \mathbf{x}^{h},t\right) :\left\Vert s^{h}\right\Vert
_{H^{2,h}\left( Q_{h,T_{1}}\right) }^{2}= \\ 
=\dsum\limits_{i,j=1}^{N}h^{2}\dint\limits_{0}^{T_{1}}\dint\limits_{0}^{1} 
\left[ \dsum\limits_{m=0}^{2}\left( \partial _{z}^{m}s^{h}\right)
^{2}+\left( \partial _{zt}^{2}s^{h}\right) ^{2}\right] \left(
x_{i},y_{j},z,t\right) dzdt \\ 
+\dsum\limits_{i,j=1}^{N}h^{2}\dint\limits_{0}^{T_{1}}\dint\limits_{0}^{1} 
\left[ \left( \partial _{t}s^{h}\right) ^{2}+\left( s^{h}\right) ^{2}\right]
\left( x_{i},y_{j},z,t\right) dzdt<\infty ,%
\end{array}%
\right.
\end{equation*}%
\begin{equation*}
H^{1,h}\left( \Omega _{h}\right) =
\end{equation*}%
\begin{equation*}
=\left\{ s^{h}\left( \mathbf{x}^{h}\right) :\left\Vert s^{h}\right\Vert
_{H^{1,h}\left( \Omega _{h}\right)
}^{2}=\dsum\limits_{i,j=1}^{N}\dsum\limits_{m=0}^{n}h^{2}\dint\limits_{0}^{1}%
\left[ \left( \partial _{z}s^{h}\right) ^{2}+\left( s^{h}\right) ^{2}\right]
\left( x_{i},y_{j},z\right) dz<\infty \right\} ,
\end{equation*}%
\begin{equation*}
L_{2}^{h}\left( \Omega _{h}\right) =\left\{ s^{h}\left( \mathbf{x}%
^{h}\right) :\left\Vert s^{h}\right\Vert _{L_{2}^{h}\left( \Omega
_{h}\right) }^{2}=\dsum\limits_{i,j=1}^{N}h^{2}\dint\limits_{0}^{1}\left(
s^{h}\right) ^{2}\left( x_{i},y_{j},z\right) dz\right\} <\infty ,
\end{equation*}%
\begin{equation*}
H^{1,h}\left( \Gamma _{h,T_{1}}\right) =\left\{ 
\begin{array}{c}
s^{h}\left( \mathbf{x}^{h},t\right) :\left\Vert s^{h}\right\Vert
_{H^{1,h}\left( \Gamma _{rh,T_{1}}\right) }^{2} \\ 
=\dsum\limits_{i,j=1}^{N}h^{2}\dint\limits_{0}^{T_{1}}\left[ \left(
s_{t}^{h}\right) ^{2}+\left( s^{h}\right) ^{2}\right] \left(
x_{i},y_{j},0,t\right) dt<\infty%
\end{array}%
\right. ,
\end{equation*}%
\begin{equation*}
L_{2}^{h}\left( \Gamma _{h,T_{1}}\right) =\left\{ s^{h}\left( \mathbf{x}%
^{h},t\right) :\left\Vert s^{h}\right\Vert _{L_{2}^{h}\left( \Gamma
_{h,T_{1}}\right)
}^{2}=\dsum\limits_{i,j=1}^{N}h\dint\limits_{0}^{T_{1}}\left( s^{h}\right)
^{2}\left( x_{i},y_{j},0,t\right) dt<\infty \right\} ,
\end{equation*}%
\begin{equation*}
L_{2}^{h}\left( \Theta _{h,T_{1}}\right) =\left\{ 
\begin{array}{c}
s^{h}\left( \mathbf{x}^{h},t\right) :\left\Vert s^{h}\right\Vert
_{L_{2}^{h}\left( \Theta _{h,T_{1}}\right) }^{2}= \\ 
=\dsum\limits_{j=0}^{N+1}h\dint\limits_{0}^{T_{1}}\dint\limits_{0}^{1}\left[
\left( s^{h}\right) ^{2}\left( x_{0},y_{j},z,t\right) +\left( s^{h}\right)
^{2}\left( x_{N+1},y_{j},z,t\right) \right] dzdt \\ 
+\dsum\limits_{i=0}^{N+1}h\dint\limits_{0}^{T_{1}}\dint\limits_{0}^{1}\left[
\left( s^{h}\right) ^{2}\left( x_{i},y_{0},z,t\right) +\left( s^{h}\right)
^{2}\left( x_{i},y_{N+1},z,t\right) \right] dzdt<\infty ,%
\end{array}%
\right.
\end{equation*}%
\begin{equation*}
C^{2,h}\left( \overline{Q_{h,T_{1}}}\right) =\left\{ 
\begin{array}{c}
s^{h}\left( \mathbf{x}^{h},t\right) :\left\Vert s^{h}\right\Vert
_{C^{2,h}\left( \overline{Q_{h,T_{1}}}\right) }= \\ 
=\max_{m=0,1,2}\left( \max_{\overline{Q_{h,T_{1}}}}\left( \left\vert
\partial _{z}^{m}s^{h}\right\vert \right) ,\max_{\overline{Q_{h,T_{1}}}%
}\left\vert \partial _{zt}^{2}s^{h}\right\vert \right) <\infty ,%
\end{array}%
\right.
\end{equation*}%
\begin{equation*}
C^{n,h}\left( \overline{\Omega _{h}}\right) =\left\{ 
\begin{array}{c}
s^{h}\left( \mathbf{x}^{h}\right) :\left\Vert s^{h}\right\Vert
_{C^{n,h}\left( \overline{\Omega _{h}}\right) } \\ 
=\max_{m\in \left[ 0,n\right] }\left( \max_{\overline{\Omega _{h}}%
}\left\vert \partial _{z}^{m}s^{h}\right\vert ,\max_{\overline{\Omega _{h}}%
}\left\vert s^{h}\right\vert \right) <\infty , \\ 
n\in \left[ 0,2\right].%
\end{array}%
\right.
\end{equation*}

\subsection{The Boundary Value Problem in Partial Finite Differences}

\label{sec:4.3}

We now rewrite the BVP (\ref{3.40})-(\ref{3.46}) as the BVP in partial
finite differences with respect to the vector functions $\left( w^{h}\left( 
\mathbf{x}^{h},t\right) ,\tau ^{h}\left( \mathbf{x}^{h}\right) \right) $:%
\begin{equation}
\Delta ^{h}w^{h}-2\tau _{z}^{h}w_{zt}^{h}-2w_{tx}^{h}\tau
_{x}^{h}-2w_{ty}^{h}\tau _{y}^{h}+  \label{4.10}
\end{equation}%
\begin{equation*}
+2w_{t}^{h}\left( \nabla ^{h}\ln w^{h}\left( \mathbf{x}^{h},0\right) \right)
\cdot \nabla ^{h}\tau ^{h}=0,\text{ }\left( \mathbf{x}^{h}\mathbf{,}t\right)
\in Q_{T_{1}}^{h},
\end{equation*}%
\begin{equation}
\Delta ^{h}\tau ^{h}+2\left( \nabla ^{h}\ln w^{h}\left( \mathbf{x}%
^{h},0\right) \right) \cdot \nabla ^{h}\tau ^{h}=0,\text{ }\mathbf{x}^{h}\in
\Omega ^{h}.  \label{4.100}
\end{equation}%
In addition, (\ref{3.45}) and (\ref{3.46}) lead to the following boundary
conditions for the $2\times 2$ system (\ref{4.10}), (\ref{4.100}) 
\begin{equation}
w^{h}\mid _{\Gamma _{h,T_{1}}}=g_{0}^{h},\text{ }w_{z}^{h}\mid _{\Gamma
_{h,T_{1}}}=g_{1}^{h},\text{ }w^{h}\mid _{\Theta _{h,T_{1}}}=g_{2}^{h},
\label{4.11}
\end{equation}%
\begin{equation}
\tau ^{h}\mid _{\Gamma _{h}}=0,~\partial _{z}\tau ^{h}\mid _{\Gamma
_{h}}=1,\tau ^{h}\mid _{\Theta _{h}}=z.  \label{4.12}
\end{equation}%
Also, using (\ref{3.41}), we impose the following condition on the function $%
w^{h}\left( \mathbf{x}^{h},0\right) :$ 
\begin{equation}
w^{h}\left( \mathbf{x}^{h},0\right) =A^{h}\left( \mathbf{x}^{h}\right) \geq
A_{0}=\frac{1}{2}\exp \left( -3n_{00}^{2}n_{0}^{2}/2\right) >0\text{ for }%
\mathbf{x}^{h}\in \overline{\Omega }^{h}.  \label{4.13}
\end{equation}

\textbf{Boundary Value Problem}$^{h}$\textbf{\ (}BVP$^{h}$). \emph{Find the
pair of functions }

$\left( w^{h},\tau _{s}^{h}\right) \in C^{2,h}\left( \overline{Q_{h,T_{1}}}%
\right) \times C^{2,h}\left( \overline{\Omega _{h}}\right) $\emph{\
satisfying conditions (\ref{4.10})-(\ref{4.12}), assuming that (\ref{4.13})
holds. }

\section{Two Carleman Estimates}

\label{sec:5}

In this section, we prove two Carleman estimates for operators written in
the above partial finite differences. Let the function $\xi \left( \mathbf{x}%
^{h}\right) \in C^{1,h}\left( \overline{\Omega _{h}}\right) $. Consider
three numbers $\xi _{0},\xi _{1},\xi _{2}$ such that $0<\xi _{0}<\xi _{1}$
and $\xi _{2}>0.$ We assume that 
\begin{equation}
0<\xi _{0}\leq \xi ^{h}\left( \mathbf{x}^{h}\right) \leq \xi _{1},\quad %
\mbox{ for all }\mathbf{x}^{h}\in \overline{\Omega _{h}},  \label{C1}
\end{equation}%
\begin{equation}
\xi _{2}=\max_{\overline{\Omega _{h}}}\left\vert \xi _{z}^{h}\left( \mathbf{x%
}^{h}\right) \right\vert .  \label{C10}
\end{equation}%
For functions $v^{h}\in H^{2,h}(Q_{h,T_{1}})$, we define the linear operator 
$L^{h}$ as: 
\begin{equation}
L^{h}v^{h}=\Delta ^{h}v^{h}-\xi ^{h}\left( \mathbf{x}^{h}\right) v_{zt}^{h},%
\text{ in }Q_{h,T_{1}},  \label{C3}
\end{equation}%
where the operator $\Delta ^{h}$ is defined in (\ref{4.80}).

\textbf{Theorem 5.1 }(the first Carleman estimate). \emph{Let the function} $%
\xi \left( \mathbf{x}^{h}\right) \in C^{1,h}\left( \overline{\Omega _{h}}%
\right) $ \emph{satisfies conditions (\ref{C1}), (\ref{C10}). Consider the
number }$\alpha _{0},$\emph{\ }%
\begin{equation}
\alpha _{0}=\alpha _{0}\left( \xi _{1}\right) =\frac{2}{3\xi _{1}}>0.
\label{C30}
\end{equation}%
\emph{\ Then there exists a sufficiently large number }$\lambda _{0}=\lambda
_{0}\left( T_{1},\xi _{0},\xi _{1},\xi _{2},h_{0},X\right) \geq 1$ \emph{%
depending only on listed\ parameters, such that for all }$\alpha \in \left(
0,\alpha _{0}\right] ,$\emph{\ all }$\lambda \geq \lambda _{0}$\emph{\ and
for all functions }$v\in H^{2,h}(Q_{h,T_{1}})$\emph{\ the following Carleman
estimate is valid:}

\begin{equation*}
\dsum\limits_{i,j=1}^{N}h^{2}\dint\limits_{0}^{T_{1}}\dint\limits_{0}^{1}%
\left( L^{h}v^{h}\right) ^{2}\left( x_{i},y_{j},z,t\right) e^{-2\lambda
\left( z+\alpha t\right) }dzdt
\end{equation*}%
\begin{equation*}
\geq C\lambda
\dsum\limits_{i,j=1}^{N}h^{2}\dint\limits_{0}^{T_{1}}\dint\limits_{0}^{1} 
\left[ \left( v_{z}^{h}\right) ^{2}+\left( v_{t}^{h}\right) ^{2}+\lambda
^{2}\left( v^{h}\right) ^{2}\right] \left( x_{i},y_{j},z,t\right)
e^{-2\lambda \left( z+\alpha t\right) }dzdt
\end{equation*}%
\begin{equation}
+C\lambda \dsum\limits_{i,j=1}^{N}h^{2}\dint\limits_{0}^{1}\left[ \left(
v_{z}^{h}\right) ^{2}+\lambda ^{2}\left( v^{h}\right) ^{2}\right] \left(
x_{i},y_{j},z,0\right) e^{-2\lambda z}dz  \label{C4}
\end{equation}%
\begin{equation*}
-C\lambda e^{-2\lambda \alpha T_{1}}\left( \left\Vert v_{z}^{h}\left( 
\mathbf{x}^{h},T_{1}\right) \right\Vert _{L_{2}^{h}\left( \Omega _{h}\right)
}^{2}+\lambda ^{2}\left\Vert v^{h}\left( \mathbf{x}^{h},T_{1}\right)
\right\Vert _{L_{2}^{h}\left( \Omega _{h}\right) }^{2}\right)
\end{equation*}%
\begin{equation*}
-C\left\Vert v^{h}\right\Vert _{L_{2}^{h}\left( \Theta _{h,T_{1}}\right)
}^{2}-C\lambda \left( \left\Vert v_{z}^{h}\right\Vert _{H^{1,h}\left( \Gamma
_{h,T_{1}}\right) }^{2}+\lambda ^{2}\left\Vert v^{h}\right\Vert
_{L_{2}^{h}\left( \Gamma _{h,T_{1}}\right) }^{2}\right) ,
\end{equation*}%
\emph{where the constant }$C=C\left( T_{1},\xi _{0},\xi _{1},\xi
_{2},h_{0},X,\alpha \right) >0$\emph{\ depends only on listed parameters. }

\textbf{Proof.} Here and below in this paper $C=C\left( T_{1},\xi _{0},\xi
_{1},\xi _{2},h_{0},X,\alpha \right) >0$ denotes different constants
depending only on listed parameters. By (\ref{4.5})-(\ref{4.80}), (\ref{C3})
and Young's inequality 
\begin{equation}
\left( L^{h}v^{h}\right) ^{2}\left( \mathbf{x}^{h},t\right) \geq \frac{1}{2}%
\left( \partial _{z}^{2}v^{h}-\xi ^{h}v_{zt}^{h}\right) ^{2}\left( \mathbf{x}%
^{h},t\right) -\left( \partial _{x}^{2}v^{h}+\partial _{y}^{2}v^{h}\right)
^{2}\left( \mathbf{x}^{h},t\right) .  \label{8.1}
\end{equation}%
Obviously%
\begin{equation*}
\dsum\limits_{i,j=1}^{N}h^{2}\dint\limits_{0}^{T_{1}}\dint\limits_{0}^{1}%
\left( \partial _{x}^{2}v^{h}+\partial _{y}^{2}v^{h}\right) ^{2}\left(
x_{i},y_{j},z,t\right) e^{-2\lambda \left( z+\alpha t\right) }dzdt
\end{equation*}%
\begin{equation*}
\leq C\left\Vert v^{h}\right\Vert _{L_{2}^{h}\left( \Theta _{h,T_{1}}\right)
}^{2}+C\dsum\limits_{i,j=1}^{N}h^{2}\dint\limits_{0}^{T_{1}}\dint%
\limits_{0}^{1}\left( v^{h}\right) ^{2}\left( x_{i},y_{j},z,t\right)
e^{-2\lambda \left( z+\alpha t\right) }dzdt.
\end{equation*}%
Hence, (\ref{8.1}) implies

\begin{equation*}
\dsum\limits_{i,j=1}^{N}h^{2}\dint\limits_{0}^{T_{1}}\dint\limits_{0}^{1}%
\left( L^{h}v^{h}\right) ^{2}\left( x_{i},y_{j},z,t\right) e^{-2\lambda
\left( z+\alpha t\right) }dzdt
\end{equation*}%
\begin{equation}
\geq \frac{1}{2}\dsum\limits_{i,j=1}^{N}h^{2}\dint\limits_{0}^{T_{1}}\dint%
\limits_{0}^{1}\left( \partial _{z}^{2}v^{h}-\left( n^{h}\right)
^{2}v_{zt}^{h}\right) ^{2}\left( x_{i},y_{j},z,t\right) e^{-2\lambda \left(
z+\alpha t\right) }dzdt  \label{8.2}
\end{equation}%
\begin{equation*}
-C\left\Vert v^{h}\right\Vert _{L_{2}^{h}\left( \Theta _{h,T_{1}}\right)
}^{2}-C\dsum\limits_{i,j=1}^{N}h^{2}\dint\limits_{0}^{T_{1}}\dint%
\limits_{0}^{1}\left( v^{h}\right) ^{2}\left( x_{i},y_{j},z,t\right)
e^{-2\lambda \left( z+\alpha t\right) }dzdt.
\end{equation*}%
Theorem 3.1 of \cite[Theorem 3.1]{Klib1d} implies that, given (\ref{C1}) and
(\ref{C30}), the following Carleman estimate holds:%
\begin{equation*}
\dsum\limits_{i,j=1}^{N}h^{2}\dint\limits_{0}^{T_{1}}\dint\limits_{0}^{1}%
\left( \partial _{z}^{2}v^{h}-\xi ^{h}v_{zt}^{h}\right) ^{2}\left(
x_{i},y_{j},z,t\right) e^{-2\lambda \left( z+\alpha t\right) }dzdt
\end{equation*}%
\begin{equation}
\geq C\lambda
\dsum\limits_{i,j=1}^{N}h^{2}\dint\limits_{0}^{T_{1}}\dint\limits_{0}^{1}%
\left( \left( v_{z}^{h}\right) ^{2}+\left( v_{t}^{h}\right) ^{2}+\lambda
^{2}\left( v^{h}\right) ^{2}\right) e^{-2\lambda \left( z+\alpha t\right)
}dzdt  \label{8.3}
\end{equation}%
\begin{equation*}
-C\lambda e^{-2\lambda \alpha T_{1}}\left( \left\Vert v_{z}^{h}\left( 
\mathbf{x}^{h},T_{1}\right) \right\Vert _{L_{2}^{h}\left( \Omega _{h}\right)
}^{2}+\lambda ^{2}\left\Vert v^{h}\left( \mathbf{x}^{h},T_{1}\right)
\right\Vert _{L_{2}^{h}\left( \Omega _{h}\right) }^{2}\right)
\end{equation*}%
\begin{equation*}
-C\lambda \left( \left\Vert v_{z}^{h}\right\Vert _{H^{1,h}\left( \Gamma
_{h,T_{1}}\right) }^{2}+\lambda ^{2}\left\Vert v^{h}\right\Vert
_{L_{2}^{h}\left( \Gamma _{h,T_{1}}\right) }^{2}\right) .
\end{equation*}%
Choosing a sufficiently large $\lambda _{0}=\lambda _{0}\left( T_{1},\xi
_{0},\xi _{1},\xi _{2},h_{0},X\right) \geq 1,$ setting $\lambda \geq \lambda
_{0}$ and combining (\ref{8.2}) and (\ref{8.3}), we obtain (\ref{C4}), which
is the target estimate of this theorem. $\square $

\textbf{Remarks 5.1}:

1. \emph{We now explain why the terms reflecting boundary conditions at }$%
\Gamma _{h,T_{1}}^{\prime }$\emph{\ are absent in the right hand side of (%
\ref{C4}) and why the condition }$\alpha \in \left( 0,\alpha _{0}\right] ,$%
\emph{\ with }$\alpha _{0}$\emph{\ as in (\ref{C30}) is imposed. The point
here is that the condition }$\alpha \in \left( 0,\alpha _{0}\right] $ \emph{%
ensures that those terms are non-negative. This follows immediately from the
combination of (\ref{C1}) with the formula (3.1) of \cite{Klib1d} as well as
with the following formulas in the proof of Theorem 3.1 of \cite{Klib1d}:
the formula (3.14) (third and fourth lines), the inequality just below
(3.14) and the formula (3.15).}

2.\emph{\ We also note that Carleman estimate (\ref{C4}) is valid for any
value }$T_{1}>0.$\emph{\ This is because of the presence of the negative
term in (\ref{C4}),}%
\begin{equation*}
-C\lambda e^{-2\lambda \alpha T_{1}}\left( \left\Vert v_{z}^{h}\left( 
\mathbf{x}^{h},T_{1}\right) \right\Vert _{L_{2}^{h}\left( \Omega _{h}\right)
}^{2}+\lambda ^{2}\left\Vert v^{h}\left( \mathbf{x}^{h},T_{1}\right)
\right\Vert _{L_{2}^{h}\left( \Omega _{h}\right) }^{2}\right) .
\end{equation*}

\textbf{Theorem 5.2} (the second Carleman estimate). \emph{Let the parameter 
}$\alpha $\emph{\ be the same as in Theorem 5.1. Then there exists a
sufficiently large number }$\lambda _{1}=\lambda _{1}\left( h_{0},X\right)
\geq 1$ \emph{such that} \emph{for all }$\lambda \geq \lambda _{1}$\emph{\
and for all functions }$v^{h}\in H^{2,h}(Q_{h,T_{1}})$\emph{\ the following
Carleman estimate is valid}%
\begin{equation*}
\dsum\limits_{i,j=1}^{N}h^{2}\dint\limits_{0}^{T_{1}}\dint\limits_{0}^{1}%
\left( \Delta ^{h}v^{h}\right) ^{2}\left( x_{i},y_{j},z,t\right)
e^{-2\lambda \left( z+\alpha t\right) }dzdt
\end{equation*}%
\begin{equation}
\geq C_{1}\lambda
\dsum\limits_{i,j=1}^{N}h^{2}\dint\limits_{0}^{T_{1}}\dint\limits_{0}^{1} 
\left[ \left( v_{z}^{h}\right) ^{2}+\lambda ^{2}\left( v^{h}\right) ^{2}%
\right] \left( x_{i},y_{j},z,t\right) e^{-2\lambda \left( z+\alpha t\right)
}dzdt  \label{C6}
\end{equation}%
\begin{equation*}
-C_{1}\left\Vert v^{h}\right\Vert _{L_{2}^{h}\left( \Theta _{h,T_{1}}\right)
}^{2}-C_{1}\lambda \left[ \left\Vert v_{z}^{h}\right\Vert _{L_{2}^{h}\left(
\Gamma _{h,T_{1}}\right) }^{2}+\lambda ^{2}\left\Vert v^{h}\right\Vert
_{L_{2}^{h}\left( \Gamma _{h,T_{1}}\right) }^{2}\right].
\end{equation*}%
\emph{Both constants }$\lambda _{1}\left( h_{0},X\right) \geq 1$\emph{\ and }%
$C_{1}=C_{1}\left( h_{0},X\right) >0$\emph{\ depend only on listed
parameters.}

\textbf{Proof.} \emph{\ }Here and below in this paper $C=C\left(
h_{0},X\right) >0$ denotes different constants depending only on listed
parameters. We obtain similarly with (\ref{8.2})%
\begin{equation*}
\dsum\limits_{i,j=1}^{N}h^{2}\dint\limits_{0}^{T_{1}}\dint\limits_{0}^{1}%
\left( \Delta ^{h}v^{h}\right) ^{2}\left( x_{i},y_{j},z,t\right)
e^{-2\lambda \left( z+\alpha t\right) }dzdt
\end{equation*}%
\begin{equation}
\geq \frac{1}{2}\dsum\limits_{i,j=1}^{N}h^{2}\dint\limits_{0}^{T_{1}}\dint%
\limits_{0}^{1}\left( v_{zz}^{h}\right) ^{2}\left( x_{i},y_{j},z,t\right)
e^{-2\lambda \left( z+\alpha t\right) }dzdt  \label{8.4}
\end{equation}%
\begin{equation*}
-C\left\Vert v^{h}\right\Vert _{L_{2}^{h}\left( \Theta _{h,T_{1}}\right)
}^{2}-C\dsum\limits_{i,j=1}^{N}h^{2}\dint\limits_{0}^{T_{1}}\dint%
\limits_{0}^{1}\left( v^{h}\right) ^{2}\left( x_{i},y_{j},z,t\right)
e^{-2\lambda \left( z+\alpha t\right) }dzdt.
\end{equation*}

It follows from \cite[Theorem 7.1]{KlibKol} that the following Carleman
estimate holds for sufficiently large $\lambda \geq \lambda _{0}\geq 1:$%
\begin{equation*}
\dsum\limits_{i,j=1}^{N}h^{2}\dint\limits_{0}^{T_{1}}\dint\limits_{0}^{1}%
\left( v_{zz}^{h}\right) ^{2}\left( x_{i},y_{j},z,t\right) e^{-2\lambda
\left( z+\alpha t\right) }dzdt
\end{equation*}%
\begin{equation}
\geq C\lambda
\dsum\limits_{i,j=1}^{N}h^{2}\dint\limits_{0}^{T_{1}}\dint\limits_{0}^{1} 
\left[ \left( v_{z}^{h}\right) ^{2}+\lambda ^{2}\left( v^{h}\right) ^{2}%
\right] \left( x_{i},y_{j},z,t\right) e^{-2\lambda \left( z+\alpha t\right)
}dzdt  \label{8.5}
\end{equation}%
\begin{equation*}
-C\lambda \left[ \left\Vert v_{z}^{h}\right\Vert _{L_{2}^{h}\left( \Gamma
_{h,T_{1}}\right) }^{2}+\lambda ^{2}\left\Vert v^{h}\right\Vert
_{L_{2}^{h}\left( \Gamma _{h,T_{1}}\right) }^{2}\right] .
\end{equation*}%
Combining (\ref{8.4}) and (\ref{8.5})$,$ we obtain (\ref{C6}), which is the
target estimate of this theorem. $\square $

\section{The stability Estimate for the CIP}

Consider the set $\mathcal{N}$ of functions $n\left( \mathbf{x}\right) $
defined as:%
\begin{equation}
\mathcal{N}=\left\{ n\left( \mathbf{x}\right) \text{ satisfies (\ref{3.1})-(%
\ref{2.3}) and }n_{z}\left( \mathbf{x}\right) \geq 0\text{ in }\mathbb{R}%
^{3}\right\} .  \label{5.4}
\end{equation}%
\textbf{Lemma 6.1}. \emph{The following inequality holds:}%
\begin{equation*}
\partial _{z}\tau \left( \mathbf{x}\right) \geq 1,\text{ }\mathbf{x}\in {%
\Omega },~\forall n\left( \mathbf{x}\right) \in \mathcal{N}.
\end{equation*}

\textbf{Proof}. It follows from (\ref{3.220})-(\ref{3.24}) that, along the
geodesic line $L\left( \mathbf{x}\right) ,$ 
\begin{equation*}
\frac{d}{ds}\left( \partial _{z}\tau \right) =\frac{n_{z}}{n},\text{ for }s>0%
\text{ }~
\end{equation*}%
and $\tau \mid _{s=0}=0.$ At $s=0$ the geodesic line $L\left( \mathbf{x}%
\right) $ intersects with the plane $\Sigma =\left\{ z=0\right\} $ defined
in (\ref{3.200}). Hence, $\partial _{x}\tau \mid _{s=0}=\partial _{y}\tau
\mid _{s=0}=0.$ By (\ref{2.3}) $n\left( \mathbf{x}\right) \mid _{\Sigma }=1.$
Hence, (\ref{3.2100}) implies that $\partial _{z}\tau \mid _{s=0}=1.$ Hence,
using (\ref{5.4}), we obtain%
\begin{equation*}
\partial _{z}\tau \left( \mathbf{x}\right) =1+\dint\limits_{L\left( \mathbf{x%
}\right) }\frac{n_{z}}{n}\left( \mathbf{\xi }\left( s\right) \right) ds\geq
1.\text{ }\square
\end{equation*}

Our CIP is an ill-posed problem. Therefore, to prove the desired stability
estimate, it is necessary to assume, in accordance with the well known
Tikhonov's concept of conditional correctness for ill-posed problems \cite{T}%
, that some \emph{a priori} known bounds are imposed on the functions $%
w^{h},\tau ^{h}$.

Thus, let $M>0$ be a positive number. Introduce the set of semi-discrete
functions $S^{h}=S^{h}\left( M,X,n_{0},n_{00}\right) $ as: \ 
\begin{equation*}
S^{h}=S^{h}\left( M,n_{0},n_{00}\right) =
\end{equation*}%
\begin{equation}
=\left\{ 
\begin{array}{c}
\left( w^{h},\tau ^{h}\right) \in C^{2,h}\left( \overline{Q_{h,T_{1}}}%
\right) \times C^{2,h}\left( \overline{\Omega _{h}}\right) : \\ 
\left\vert \nabla ^{h}\tau ^{h}\left( \mathbf{x}^{h}\right) \right\vert \leq
n_{0},\text{ }\mathbf{x}^{h}\in \overline{\Omega _{h}}, \\ 
\left\Vert w^{h}\right\Vert _{C^{2,h}\left( \overline{Q_{h,T_{1}}}\right)
},\left\Vert \tau ^{h}\right\Vert _{C^{2,h}\left( \overline{\Omega _{h}}%
\right) }\leq M, \\ 
1\leq \partial _{z}\tau ^{h}\left( \mathbf{x}\right) \leq n_{0},\text{ }%
\mathbf{x}^{h}\in \overline{\Omega }_{h}, \\ 
w^{h}\left( \mathbf{x}^{h},0\right) \geq A_{0}, \\ 
\text{where the vector function }\left( w^{h},\tau ^{h}\right) \\ 
\text{ is generated by the above procedure} \\ 
\text{ and the number }A_{0}\text{ is defined in (\ref{3.262}).}%
\end{array}%
\right.  \label{5.9}
\end{equation}%
Conditions in the second, fourth and fifth lines of (\ref{5.9}) are imposed
due to (\ref{2.2}), (\ref{3.21}), Lemma 6.1 and Theorem 3.1, respectively.
We have 
\begin{equation}
S^{h}=S_{w}^{h}\times S_{\tau }^{h}=\left\{ \left( w^{h}\left( \mathbf{x}%
^{h},t\right) ,\tau ^{h}\left( \mathbf{x}^{h}\right) \right) \right\} .
\label{5.90}
\end{equation}%
Thus, by (\ref{5.9}) and (\ref{5.90}) 
\begin{equation}
w^{h}\left( \mathbf{x}^{h},0\right) \geq A_{0},\mathbf{x}^{h}\in \overline{%
\Omega }_{h},\text{ }\forall w^{h}\in S_{w}^{h}.  \label{5.10}
\end{equation}%
\begin{equation}
\text{ }1\leq \partial _{z}\tau ^{h}\left( \mathbf{x}^{h}\right) \leq n_{0},%
\text{ }\mathbf{x}^{h}\in \overline{\Omega }_{h},\forall \tau ^{h}\in
S_{\tau }^{h}.  \label{5.11}
\end{equation}%
Using (\ref{3.21}), define the function $n^{h}\left( \mathbf{x}^{h}\right) $
as: 
\begin{equation}
\left\vert \nabla _{\mathbf{x}}\tau ^{h}\left( \mathbf{x}^{h}\right)
\right\vert ^{2}=\left( n^{h}\right) ^{2}\left( \mathbf{x}^{h}\right) ,\text{
}\forall \tau ^{h}\in S_{\tau }^{h}.  \label{5.80}
\end{equation}

\textbf{Theorem 6.1} (H\"{o}lder stability estimate). \emph{Let two vector
functions }$\left( w_{1}^{h},\tau _{1}^{h}\right) $, $\left( w_{2}^{h},,\tau
_{2}^{h}\right) \in S^{h}$\emph{\ be solutions of BVP}$^{h}$\emph{\ with two
sets of boundary data at }$\Gamma _{h,T_{1}},\Gamma _{h},\Theta _{h,T_{1}}$%
\emph{and} $\Theta _{h}$, 
\begin{equation}
w_{j}^{h}\mid _{\Gamma _{h,T_{1}}}=g_{j,0}^{h},\text{ }\partial
_{z}w_{j}^{h}\mid _{\Gamma _{h,T_{1}}}=g_{j,1}^{h},~w_{j}^{h}\mid _{\Theta
_{h,T_{1}}}=g_{j,2}^{h},\text{ }j=1,2,  \label{5.12}
\end{equation}%
\begin{equation}
\tau _{j}^{h}\mid _{\Gamma _{h}}=0,~\partial _{z}\tau _{j}^{h}\mid _{\Gamma
_{h}}=1,~\tau _{j}^{h}\mid _{\Theta _{h}}=z,\text{ }j=1,2.  \label{5.13}
\end{equation}

\emph{Denote}%
\begin{equation}
\left\{ 
\begin{array}{c}
\widetilde{w}^{h}=w_{1}^{h}-w_{2}^{h},~\widetilde{\tau }^{h}=\tau
_{1}^{h}-\tau _{2}^{h}, \\ 
\widetilde{g}_{0}^{h}=g_{1,0}^{h}-g_{2,0}^{h},~\widetilde{g}%
_{1}^{h}=g_{1,1}^{h}-g_{2,1}^{h},~\widetilde{g}%
_{2}^{h}=g_{1,2}^{h}-g_{2,2}^{h}, \\ 
\widetilde{n}^{h}=n_{1}^{h}-n_{2}^{h},%
\end{array}%
\right.  \label{5.14}
\end{equation}%
\emph{\ where functions }$n_{1}^{h}$\emph{\ and }$n_{2}^{h}$\emph{\ are
obtained from functions }$\tau _{1}^{h}$\emph{\ and }$\tau _{2}^{h}$\emph{\
respectively via (\ref{5.80}). Assume that}%
\begin{equation}
\left\Vert \widetilde{g}_{0}^{h}\right\Vert _{H^{1,h}\left( \Gamma
_{h,T_{1}}\right) },~\left\Vert \widetilde{g}_{1}^{h}\right\Vert
_{L_{2}^{h}\left( \Gamma _{h,T_{1}}\right) },\left\Vert \widetilde{g_{2}}%
^{h}\right\Vert _{L_{2}^{h}\left( \Theta _{h,T_{1}}\right) }<\delta ,
\label{5.140}
\end{equation}%
\emph{\ where }$\delta \in \left( 0,1\right) $\emph{\ is a number. Define
the number }$\alpha _{0}=\alpha _{0}\left( n_{0}\right) =2/\left(
3n_{0}\right) .$ \emph{Consider an arbitrary number} $\alpha \in \left(
0,\alpha _{0}\right] $\emph{\ and set }$T_{1}=3/\alpha $. \emph{Then there
exists a sufficiently small number }$\delta _{0}=\delta _{0}\left(
M,n_{0},n_{00},h_{0},X,\alpha \right) \in \left( 0,1\right) $\emph{\ and a
number }$C_{2}=C_{2}\left( M,n_{0},n_{00},h_{0},X,\alpha \right) >0,$\emph{\
both numbers depending only on listed parameters, such that if }$\delta \in
\left( 0,\delta _{0}\right) ,$\emph{\ then the following stability estimates
are valid for the functions} $\widetilde{w}^{h},\widetilde{\tau }^{h},%
\widetilde{n}^{h}:$%
\begin{equation}
\left\Vert \widetilde{w}^{h}\right\Vert _{H^{1,h}\left( Q_{h,1/\alpha
}\right) },~\left\Vert \widetilde{\tau }^{h}\right\Vert _{H^{1,h}\left(
\Omega _{h}\right) },~\left\Vert \widetilde{n}^{h}\right\Vert
_{L_{2}^{h}\left( \Omega _{h}\right) }\leq C_{2}\delta^{1/3}\ln\left( \delta
^{-1}\right).  \label{5.17}
\end{equation}

\textbf{Proof.} Below $C_{2}=C_{2}\left( M,n_{0},n_{00},h_{0},X,\alpha
\right) >0$ denotes different constants depending only on listed parameters.
Let $a_{1},a_{2},b_{1},b_{2}\in \mathbb{R}$ be arbitrary numbers. Let $%
\widetilde{a}=a_{1}-a_{2}$ and $\widetilde{b}=b_{1}-b_{2}.$ Then \ 
\begin{equation}
a_{1}b_{1}-a_{2}b_{2}=\widetilde{a}b_{1}+\widetilde{b}a_{2}.  \label{6.0}
\end{equation}%
Using (\ref{3.262}), (\ref{5.9}), (\ref{5.14}) and (\ref{6.0}), we obtain 
\begin{equation}
\left\vert \nabla ^{h}\ln w_{1}^{h}\left( \mathbf{x}^{h},0\right) -\nabla
^{h}\ln w_{2}^{h}\left( \mathbf{x}^{h},0\right) \right\vert \leq C_{2}\left(
\left\vert \nabla ^{h}\widetilde{w}\left( \mathbf{x}^{h},0\right)
\right\vert +\left\vert \widetilde{w}\left( \mathbf{x}^{h},0\right)
\right\vert \right) .  \label{6.00}
\end{equation}

It is well known that, when applying a Carleman estimate, one can replace
differential equations with appropriate differential inequalities, see, e.g. 
\cite{KL}. Hence, subtract two equations (\ref{4.10}), (\ref{4.100}) for the
pair $\left( w_{2}^{h},,\tau _{2}^{h}\right) $ from two equations (\ref{4.10}%
), (\ref{4.100}) for the pair $\left( w_{1}^{h},\tau _{1}^{h}\right) ,$ use
the first line of (\ref{5.14}) and (\ref{6.0}). Then, using (\ref{6.00}),
turn resulting equations in inequalities with respect to functions $%
\widetilde{w}^{h}$ and $\widetilde{\tau }_{s}^{h}.$ We obtain two
differential inequalities: 
\begin{equation*}
\left\vert \Delta ^{h}\widetilde{w}^{h}-\partial _{z}\tau _{1}^{h}\widetilde{%
w}_{zt}^{h}\right\vert
\end{equation*}%
\begin{equation}
\leq C_{2}\left( \left\vert \widetilde{w}_{t}^{h}\right\vert +\left\vert
\nabla ^{h}\widetilde{w}^{h}\left( \mathbf{x}^{h}\mathbf{,}0\right)
\right\vert +\left\vert \widetilde{w}^{h}\left( \mathbf{x}^{h}\mathbf{,}%
0\right) \right\vert +\left\vert \nabla ^{h}\widetilde{\tau }^{h}\right\vert
\right) ,\text{ }\left( \mathbf{x}^{h}\mathbf{,}t\right) \in Q_{T_{1}}^{h},
\label{6.1}
\end{equation}%
\begin{equation}
\left\vert \Delta ^{h}\widetilde{\tau }^{h}\right\vert \leq C_{2}\left(
\left\vert \nabla ^{h}\widetilde{w}^{h}\left( \mathbf{x}^{h}\mathbf{,}%
0\right) \right\vert +\left\vert \widetilde{w}^{h}\left( \mathbf{x}^{h}%
\mathbf{,}0\right) \right\vert +\left\vert \nabla ^{h}\widetilde{\tau }%
^{h}\right\vert \right) ,\text{ }\mathbf{x}^{h}\in \Omega _{h}.  \label{6.2}
\end{equation}%
By (\ref{5.13}) and the second line of (\ref{5.14}) the boundary conditions
are:%
\begin{equation}
\widetilde{w}^{h}\mid _{\Gamma _{h,T_{1}}}=\widetilde{g}_{0}^{h},\text{ }%
\widetilde{w}_{z}^{h}\mid _{\Gamma _{h,T_{1}}}=\widetilde{g}_{1}^{h},\text{ }%
\widetilde{w}^{h}\mid _{\Theta _{h,T_{1}}}=\widetilde{g}_{2}^{h},
\label{6.3}
\end{equation}%
\begin{equation}
\widetilde{\tau }^{h}\mid _{\Gamma _{h}}=0,~\partial _{z}\widetilde{\tau }%
^{h}\mid _{\Gamma _{h}}=0,~\widetilde{\tau }^{h}\mid _{\Theta _{h}}=0.
\label{6.4}
\end{equation}

Denote%
\begin{equation*}
I\left( x_{i},y_{j},z,t\right) =
\end{equation*}%
\begin{equation}
\left[ \left\vert \widetilde{w}_{t}^{h}\left( x_{i},y_{j},z,t\right)
\right\vert +\left( \left\vert \nabla ^{h}\widetilde{w}^{h}\right\vert
+\left\vert \widetilde{w}^{h}\right\vert \right) \left(
x_{i},y_{j},z,0\right) +\left\vert \nabla ^{h}\widetilde{\tau }\right\vert
\left( x_{i},y_{j},z\right) \right] ^{2}.  \label{6.40}
\end{equation}%
Using Cauchy-Schwarz inequality, we obtain%
\begin{equation*}
\dsum\limits_{i,j=1}^{N}h^{2}\dint\limits_{0}^{T_{1}}\dint\limits_{0}^{1}I%
\left( x_{i},y_{j},z,t\right) e^{-2\lambda \left( z+\alpha t\right) }dzdt
\end{equation*}%
\begin{equation*}
\leq
C_{1}\dsum\limits_{i,j=1}^{N}h^{2}\dint\limits_{0}^{T_{1}}\dint%
\limits_{0}^{1}\left( \widetilde{w}_{t}^{h}\left( x_{i},y_{j},z,t\right)
\right) ^{2}e^{-2\lambda \left( z+\alpha t\right) }dzdt
\end{equation*}%
\begin{equation}
+C_{1}\dsum\limits_{i,j=1}^{N}h^{2}\dint\limits_{0}^{T_{1}}\dint%
\limits_{0}^{1}\left[ \left( \left( w_{z}^{h}\right) ^{2}+\left(
w^{h}\right) ^{2}\right) \left( x_{i},y_{j},z,0\right) \right] e^{-2\lambda
\left( z+\alpha t\right) }dzdt  \label{6.41}
\end{equation}%
\begin{equation*}
+C_{1}\dsum\limits_{i,j=1}^{N}h^{2}\dint\limits_{0}^{T_{1}}\dint%
\limits_{0}^{1}\left[ \left( \partial _{z}\widetilde{\tau }^{h}\right)
^{2}+\left( \widetilde{\tau }^{h}\right) ^{2}\right] \left(
x_{i},y_{j},z\right) e^{-2\lambda \left( z+\alpha t\right) }dzdt.
\end{equation*}%
And also%
\begin{equation*}
\dsum\limits_{i,j=1}^{N}h^{2}\dint\limits_{0}^{T_{1}}\dint\limits_{0}^{1} 
\left[ \left( \left\vert \nabla ^{h}\widetilde{w}^{h}\right\vert +\left\vert 
\widetilde{w}^{h}\right\vert \right) \left( x_{i},y_{j},z,0\right)
+\left\vert \nabla ^{h}\widetilde{\tau }\right\vert \left(
x_{i},y_{j},z\right) \right] ^{2}e^{-2\lambda \left( z+\alpha t\right) }dzdt
\end{equation*}%
\begin{equation}
\leq
C_{1}\dsum\limits_{i,j=1}^{N}h^{2}\dint\limits_{0}^{T_{1}}\dint%
\limits_{0}^{1}\left[ \left( \partial _{z}\widetilde{\tau }^{h}\right)
^{2}+\left( \widetilde{\tau }^{h}\right) ^{2}\right] \left(
x_{i},y_{j},z\right) e^{-2\lambda \left( z+\alpha t\right) }dzdt
\label{6.42}
\end{equation}%
\begin{equation*}
+C_{1}\dsum\limits_{i,j=1}^{N}h^{2}\dint\limits_{0}^{T_{1}}\dint%
\limits_{0}^{1}\left[ \left( \left( w_{z}^{h}\right) ^{2}+\left(
w^{h}\right) ^{2}\right) \left( x_{i},y_{j},z,0\right) \right] e^{-2\lambda
\left( z+\alpha t\right) }dzdt.
\end{equation*}

Square both sides of both inequalities (\ref{6.1}) and (\ref{6.2}). Then
multiply the results by the function $e^{-2\lambda \left( z+\alpha t\right)
} $, construct sums combined with integrals like in the left hand sides of
Carleman estimates (\ref{C4}) and (\ref{C6}) and then use (\ref{6.40})-(\ref%
{6.42}). We obtain 
\begin{equation*}
\dsum\limits_{i,j=1}^{N}h^{2}\dint\limits_{0}^{T_{1}}\dint\limits_{0}^{1}%
\left( \Delta ^{h}\widetilde{w}^{h}-\partial _{z}\tau _{1}^{h}\widetilde{w}%
_{zt}^{h}\right) ^{2}\left( x_{i},y_{j},z,t\right) e^{-2\lambda \left(
z+\alpha t\right) }dzdt
\end{equation*}%
\begin{equation*}
\leq
C_{1}\dsum\limits_{i,j=1}^{N}h^{2}\dint\limits_{0}^{T_{1}}\dint%
\limits_{0}^{1}\left( \widetilde{w}_{t}^{h}\left( x_{i},y_{j},z,t\right)
\right) ^{2}e^{-2\lambda \left( z+\alpha t\right) }dzdt
\end{equation*}%
\begin{equation}
+C_{1}\dsum\limits_{i,j=1}^{N}h^{2}\dint\limits_{0}^{T_{1}}\dint%
\limits_{0}^{1}\left[ \left( \widetilde{w}_{z}^{h}\right) ^{2}+\left( 
\widetilde{w}^{h}\right) ^{2}\right] \left( x_{i},y_{j},z,0\right)
e^{-2\lambda \left( z+\alpha t\right) }dzdt  \label{6.5}
\end{equation}%
\begin{equation*}
+C_{1}\dsum\limits_{i,j=1}^{N}h^{2}\dint\limits_{0}^{T_{1}}\dint%
\limits_{0}^{1}\left[ \left( \partial _{z}\widetilde{\tau }^{h}\right)
^{2}+\left( \widetilde{\tau }^{h}\right) ^{2}\right] \left(
x_{i},y_{j},z\right) e^{-2\lambda \left( z+\alpha _{0}t\right) }dzdt.
\end{equation*}%
And also 
\begin{equation*}
\dsum\limits_{i,j=1}^{N}h^{2}\dint\limits_{0}^{T_{1}}\dint\limits_{0}^{1}%
\left( \Delta ^{h}\widetilde{\tau }^{h}\right) ^{2}\left(
x_{i},y_{j},z\right) e^{-2\lambda \left( z+\alpha t\right) }dzdt
\end{equation*}%
\begin{equation}
\leq
C_{1}\dsum\limits_{i,j=1}^{N}h^{2}\dint\limits_{0}^{T_{1}}\dint%
\limits_{0}^{1}\left[ \left( \partial _{z}\widetilde{\tau }^{h}\right)
^{2}+\left( \widetilde{\tau }^{h}\right) ^{2}\right] \left(
x_{i},y_{j},z\right) e^{-2\lambda \left( z+\alpha t\right) }dzdt  \label{6.6}
\end{equation}%
\begin{equation*}
+C_{1}\dsum\limits_{i,j=1}^{N}h^{2}\dint\limits_{0}^{T_{1}}\dint%
\limits_{0}^{1}\left[ \left( \widetilde{w}_{z}^{h}\right) ^{2}+\left( 
\widetilde{w}^{h}\right) ^{2}\right] \left( x_{i},y_{j},z,0\right)
e^{-2\lambda \left( z+\alpha t\right) }dzdt.
\end{equation*}
By (\ref{C30}) and the second line of (\ref{5.9}) Set in Theorem 5.1 $\xi
^{h}\left( \mathbf{x}^{h}\right) =\partial _{z}\tau _{1}^{h}\left( \mathbf{x}%
^{h}\right) .$ Then Theorem 5.1 is applicable here. Indeed, it follows from (%
\ref{5.9}) and Lemma 6.1 that we can now take in (\ref{C1}) $\xi _{0}=1$, $%
\xi _{1}=n_{0},$ and we can take $\xi _{2}=M$ in (\ref{C10}). Let $\lambda
_{0}\geq 1$\ be the same as in Theorem 5.1. Thus, apply Carleman estimates (%
\ref{C4}) and (\ref{C6}) to the left hand sides of (\ref{6.5}) and (\ref{6.6}%
) respectively. Then sum up resulting inequalities. In doing so, we replace
in (\ref{C4}) $v^{h}\left( \mathbf{x}^{h},t\right) $ with $\widetilde{w}%
^{h}\left( \mathbf{x}^{h},t\right) $ and replace in (\ref{C6}) $v^{h}\left( 
\mathbf{x}^{h},t\right) $ with $\widetilde{\tau }^{h}\left( \mathbf{x}%
^{h}\right) .$ In addition, use (\ref{6.4}). We obtain for all $\lambda \geq
\lambda _{0}$%
\begin{equation*}
\dsum\limits_{i,j=1}^{N}h^{2}\dint\limits_{0}^{T_{1}}\dint\limits_{0}^{1}%
\left( \widetilde{w}_{t}^{h}\left( x_{i},y_{j},z,t\right) \right)
^{2}e^{-2\lambda \left( z+\alpha t\right) }dzdt
\end{equation*}%
\begin{equation*}
+\dsum\limits_{i,j=1}^{N}h^{2}\dint\limits_{0}^{T_{1}}\dint\limits_{0}^{1} 
\left[ \left( \widetilde{w}^{h}\right) ^{2}+\left( \widetilde{w}%
_{z}^{h}\right) ^{2}\right] \left( x_{i},y_{j},z,0\right) e^{-2\lambda
\left( z+\alpha t\right) }dzdt
\end{equation*}%
\begin{equation*}
+\dsum\limits_{i,j=1}^{N}h^{2}\dint\limits_{0}^{T_{1}}\dint\limits_{0}^{1} 
\left[ \left( \partial _{z}\widetilde{\tau }^{h}\right) ^{2}+\left( 
\widetilde{\tau }^{h}\right) ^{2}\right] \left( x_{i},y_{j},z\right)
e^{-2\lambda \left( z+\alpha t\right) }dzdt.
\end{equation*}%
\begin{equation*}
\geq C_{2}\lambda
\dsum\limits_{i,j=1}^{N}h^{2}\dint\limits_{0}^{T_{1}}\dint\limits_{0}^{1}%
\left( \left( \widetilde{w}_{z}^{h}\right) ^{2}+\left( \widetilde{w}%
_{t}^{h}\right) ^{2}+\lambda ^{2}\left( \widetilde{w}^{h}\right) ^{2}\right)
\left( x_{i},y_{j},z,t\right) e^{-2\lambda \left( z+\alpha t\right) }dzdt
\end{equation*}%
\begin{equation}
+C_{2}\lambda
\dsum\limits_{i,j=1}^{N}h^{2}\dint\limits_{0}^{T_{1}}\dint\limits_{0}^{1}%
\left( \left( \partial _{z}\widetilde{\tau }^{h}\right) ^{2}+\lambda
^{2}\left( \widetilde{\tau }^{h}\right) ^{2}\right) \left(
x_{i},y_{j},z\right) e^{-2\lambda \left( z+\alpha t\right) }dzdt  \label{6.7}
\end{equation}%
\begin{equation*}
+C_{2}\lambda \dsum\limits_{i,j=1}^{N}h^{2}\dint\limits_{0}^{1}\left[ \left( 
\widetilde{w}_{z}^{h}\right) ^{2}+\lambda ^{2}\left( \widetilde{w}%
^{h}\right) ^{2}\right] \left( x_{i},y_{j},z,0\right) e^{-2\lambda z}dz
\end{equation*}%
\begin{equation*}
-C_{2}\lambda ^{3}e^{-6\lambda }\left\Vert \widetilde{w}^{h}\left( \mathbf{x}%
^{h},T_{1}\right) \right\Vert _{H^{1,h}\left( \Omega _{h}\right)
}^{2}-C_{2}\lambda ^{3}\left( \left\Vert \widetilde{w}^{h}\right\Vert
_{H^{1,h}\left( \Gamma _{h,T_{1}}\right) }^{2}+\left\Vert \widetilde{w}%
_{z}^{h}\right\Vert _{L_{2}^{h}\left( \Gamma _{h,T1}\right) }^{2}\right)
\end{equation*}%
\begin{equation*}
-C_{2}\left\Vert \widetilde{w}^{h}\right\Vert _{L_{2}\left( \Theta
_{h,T_{1}}\right) }^{2}.
\end{equation*}%
The multiplier $e^{-6\lambda }$ in the sixth line of (\ref{6.7}) is due to
the fact that $T_{1}=3/\alpha $ and $e^{-2\lambda \left( z+\alpha
T_{1}\right) }\leq e^{-6\lambda }$ for $z\in \left( 0,1\right) .$ Choose $%
\lambda _{1}=\lambda _{1}\left( M,n_{0},n_{00},h_{0},X,\alpha \right) \geq
\lambda _{0}\geq 1$ so large that $C_{2}\lambda _{1}/2>1.$ Then (\ref{5.14}%
), (\ref{6.3}) and (\ref{6.7}) imply%
\begin{equation*}
C_{2}\lambda ^{3}e^{-6\lambda }\left\Vert \widetilde{w}^{h}\left( \mathbf{x}%
^{h},T_{1}\right) \right\Vert _{H^{1,h}\left( \Omega _{h}\right)
}^{2}+C_{2}\lambda ^{3}\left( \left\Vert \widetilde{g}_{0}^{h}\right\Vert
_{H^{1,h}\left( \Gamma _{h,T_{1}}\right) }^{2}+\left\Vert \widetilde{g}%
_{1}^{h}\right\Vert _{L_{2}^{h}\left( \Gamma _{h,T1}\right) }^{2}\right)
\end{equation*}%
\begin{equation*}
+C_{2}\left( \left\Vert \widetilde{w}^{h}\right\Vert _{L_{2}\left( \Theta
_{h,T_{1}}\right) }^{2}+\left\Vert \widetilde{w}^{h}\right\Vert
_{L_{2}\left( \Theta _{h}\right) }^{2}\right)
\end{equation*}%
\begin{equation}
\geq \lambda
\dsum\limits_{i,j=1}^{N}h^{2}\dint\limits_{0}^{T_{1}}\dint\limits_{0}^{1}%
\left( \left( \widetilde{w}_{z}^{h}\right) ^{2}+\left( \widetilde{w}%
_{t}^{h}\right) ^{2}+\lambda ^{2}\left( \widetilde{w}^{h}\right) ^{2}\right)
\left( x_{i},y_{j},z,t\right) e^{-2\lambda \left( z+\alpha _{0}t\right) }dzdt
\label{6.70}
\end{equation}%
\begin{equation*}
+\lambda
\dsum\limits_{i,j=1}^{N}h^{2}\dint\limits_{0}^{T_{1}}\dint\limits_{0}^{1}%
\left( \left( \partial _{z}\widetilde{\tau }^{h}\right) ^{2}+\lambda
^{2}\left( \widetilde{\tau }^{h}\right) ^{2}\right) \left(
x_{i},y_{j},z\right) e^{-2\lambda \left( z+\alpha t\right) }dzdt
\end{equation*}%
\begin{equation*}
+\lambda \dsum\limits_{i,j=1}^{N}h^{2}\dint\limits_{0}^{1}\left[ \left( 
\widetilde{w}_{z}^{h}\right) ^{2}+\lambda ^{2}\left( \widetilde{w}%
^{h}\right) ^{2}\right] \left( x_{i},y_{j},z,0\right) e^{-2\lambda z}dz,%
\text{ }\forall \lambda \geq \lambda _{1}.
\end{equation*}

Introduce the number $t_{1}=1/\alpha \in \left( 0,T_{1}\right) =\left(
0,3/\alpha \right) .$ Then 
\begin{equation}
e^{-2\lambda \left( z+\alpha t\right) }>e^{-2\lambda \left( 1+\alpha
t_{1}\right) }=e^{-4\lambda }\text{ for }\left( z,t\right) \in \left(
0,1\right) \times \left( 0,t_{1}\right) .  \label{6.71}
\end{equation}%
Divide (\ref{6.70}) by $\lambda e^{-4\lambda }$ and replace $\lambda
^{2}\geq 1$ with $1$ in the right hand side of the resulting inequality,
thus, making it stronger. Using (\ref{6.71}), we obtain 
\begin{equation*}
C_{2}\lambda ^{2}e^{4\lambda }\left( \left\Vert \widetilde{g}%
_{0}^{h}\right\Vert _{H^{1,h}\left( \Gamma _{h,T_{1}}\right)
}^{2}+\left\Vert \widetilde{g}_{0}^{h}\right\Vert _{L_{2}^{h}\left( \Gamma
_{h,T_{1}}\right) }^{2}\right)
\end{equation*}%
\begin{equation*}
+C_{2}e^{4\lambda }\left( \left\Vert \widetilde{w}^{h}\right\Vert
_{L_{2}\left( \Theta _{h,T_{1}}\right) }^{2}+\left\Vert \widetilde{w}%
^{h}\right\Vert _{L_{2}\left( \Theta _{h}\right) }^{2}\right) +C_{2}\lambda
^{2}e^{-2\lambda }\left\Vert \widetilde{w}^{h}\left( \mathbf{x}%
^{h},T_{1}\right) \right\Vert _{H^{1,h}\left( \Omega _{h}\right) }^{2}
\end{equation*}%
\begin{equation}
\geq
\dsum\limits_{i,j=1}^{N}h^{2}\dint\limits_{0}^{t_{1}}\dint\limits_{0}^{1}%
\left( \left( \widetilde{w}_{z}^{h}\right) ^{2}+\left( \widetilde{w}%
_{t}^{h}\right) ^{2}+\left( \widetilde{w}^{h}\right) ^{2}\right) \left(
x_{i},y_{j},z,t\right) dzdt  \label{6.8}
\end{equation}%
\begin{equation*}
+\dsum\limits_{i,j=1}^{N}h^{2}\dint\limits_{0}^{t_{1}}\dint\limits_{0}^{1}%
\left( \left( \partial _{z}\widetilde{\tau }^{h}\right) ^{2}+\left( 
\widetilde{\tau }^{h}\right) ^{2}\right) \left( x_{i},y_{j},z\right) dzdt
\end{equation*}%
\begin{equation*}
+\dsum\limits_{i,j=1}^{N}h^{2}\dint\limits_{0}^{1}\left[ \left( \widetilde{w}%
_{z}^{h}\right) ^{2}+\lambda ^{2}\left( \widetilde{w}^{h}\right) ^{2}\right]
\left( x_{i},y_{j},z,0\right) e^{-2\lambda z}dz,\text{ }\forall \lambda \geq
\lambda _{1}.
\end{equation*}%
Choose $\delta _{0}=\delta _{0}\left( M,n_{0},n_{00},h_{0},X,\alpha \right)
\in \left( 0,1\right) $ such that $\lambda _{1}=\ln\left(\delta
_{0}^{-1/3}\right).$ For every $\delta \in \left( 0,\delta _{0}\right) $ we
choose $\lambda =\lambda \left( \delta \right) >\lambda _{1}$ such that 
\begin{equation}
e^{4\lambda }\delta ^{2}=e^{-2\lambda }.  \label{6.90}
\end{equation}

Hence, $\lambda =\lambda \left( \delta \right) =\frac{1}{3}\ln \left( \delta
^{-1}\right) .$ It follows from the third line of (\ref{5.9}) and the first
line of (\ref{5.14}) that $\left\Vert \widetilde{w}^{h}\left( \mathbf{x}%
^{h},T_{1}\right) \right\Vert _{H^{1,h}\left( \Omega _{h}\right) }^{2}\leq
C_{2}$. Hence, using (\ref{5.140}), (\ref{6.8}) and (\ref{6.90}), we obtain 
\begin{equation}
\left\Vert \widetilde{w}^{h}\right\Vert _{H^{1,h}\left( Q_{h,1/\alpha
}\right) },\left\Vert \widetilde{\tau }\right\Vert _{H^{1,h}\left( \Omega
_{h}\right) }\leq C_{2}\delta^{1/3}\ln\left( \delta ^{-1}\right).
\label{6.100}
\end{equation}%
Next, by the fourth line of (\ref{5.9}) and (\ref{5.14}) 
\begin{equation*}
\left\vert \widetilde{n}^{h}\right\vert =\left\vert \frac{%
(n_{1}^{h})^{2}-(n_{2}^{h})^{2}}{n_{1}^{h}+n_{2}^{h}}\right\vert \leq \frac{1%
}{2}\left\vert (n_{1}^{h})^{2}-(n_{2}^{h})^{2}\right\vert =\frac{1}{2}%
\left\vert |\nabla \tau _{1}^{h}|^{2}-|\nabla \tau _{1}^{h}|^{2}\right\vert
\end{equation*}%
\begin{equation*}
\leq \frac{1}{2}\left\vert \nabla ^{h}\widetilde{\tau }\cdot (\nabla \tau
_{1}^{h}+\nabla \tau _{2}^{h})\right\vert \leq n_{0}\left\vert \nabla ^{h}%
\widetilde{\tau }\right\vert.
\end{equation*}%
Hence, by (\ref{6.100}) 
\begin{equation}
\left\Vert \widetilde{n}^{h}\right\Vert _{L_{2}^{h}\left( \Omega _{h}\right)
}\leq C_{2}\delta^{1/3}\ln\left( \delta ^{-1}\right).  \label{6.110}
\end{equation}%
Estimates (\ref{6.100}) and (\ref{6.110}) imply the target estimates (\ref%
{5.17}). $\ \square $

\begin{center}
\textbf{Acknowledgment}
\end{center}

The work of V.G. Romanov was supported by a grant from the Siberian Branch
of the Russian Academy of Science, project number FWNF-2022-0009.

\end{document}